\documentclass{amsart}

\usepackage{amssymb,latexsym,amsmath,amsfonts,hhline,latexsym}
\usepackage{pdfsync}

\newcommand{\cal}{\mathcal}
\usepackage{cite}
\usepackage[matrix,arrow,curve]{xy}
\input xy
\xyoption{all}

\oddsidemargin=\evensidemargin \advance\topmargin by-\headheight

\def\VRT#1{*=<7mm>[o][F-]{#1}}

\DeclareMathOperator{\aut}{Aut}

\DeclareMathOperator{\cyc}{Cyc}

\DeclareMathOperator{\Ext}{Ext}

\DeclareMathOperator{\id}{id}
\DeclareMathOperator{\hol}{Hol}

\DeclareMathOperator{\iso}{Iso}

\DeclareMathOperator{\mult}{Mult}

\DeclareMathOperator{\orb}{Orb}

\DeclareMathOperator{\rad}{rad}
\DeclareMathOperator{\rk}{rk}

\DeclareMathOperator{\sch}{Sch}
\DeclareMathOperator{\Span}{span}

\DeclareMathOperator{\sym}{Sym}

\def\mZ{{\mathbb Z}}

\def\frS{{\mathfrak S}}
\def\frSH{{\widehat\frS}}
\def\fS{{\Sigma}}

\def\cA{{\cal A}}
\def\cB{{\cal B}}

\def\cF{{\frS}}
\def\cG{{\cal G}}
\def\cH{{\cal G}}

\def\cK{{\cal K}}
\def\cM{{\cal M}}
\def\cP{{\cal P}}

\def\cS{{\cal S}}

\def\lg{\langle}

\def\rg{\rangle}
\def\twoe{\underset{\scriptscriptstyle ^2}{\approx}}
\def\wh{\widehat}

\def\proof{{\bf Proof}.\ }
\def\bull{\vrule height .9ex width .8ex depth -.1ex }
\def\qaq{\quad\text{and}\quad}
\newcommand{\grp}[1]{\langle {#1}\rangle}

\makeatletter
\renewcommand{\subsection}{\@startsection{subsection}{2}{0mm}{-2mm}{-2mm}{\bf\normalsize}}

\def\sbsnt#1{\subsection{#1}}
\makeatother

\newtheorem{formula}{}[section]
\newtheorem{proposition}[formula]{Proposition}
\newtheorem{definition}[formula]{Definition}
\newtheorem{corollary}[formula]{Corollary}
\newtheorem{remark}[formula]{Remark}
\newtheorem{lemma}[formula]{Lemma}
\newtheorem{theorem}[formula]{Theorem}

\newtheorem{example}[formula]{Example}

\def\thrm{\begin{theorem}}
\def\thrml#1{\begin{theorem}\label{#1}}
\def\ethrm{\end{theorem}}
\def\prpstn{\begin{proposition}}
\def\prpstnl#1{\begin{proposition}\label{#1}}
\def\eprpstn{\end{proposition}}
\def\rmrk{\begin{remark}}
\def\rmrkl#1{\begin{remark}\label{#1}}
\def\ermrk{\end{remark}}
\def\dfntn{\begin{definition}}
\def\dfntnl#1{\begin{definition}\label{#1}}
\def\edfntn{\end{definition}}
\def\nmrt{\begin{enumerate}}
\def\enmrt{\end{enumerate}}
\def\tm#1{\item[{\rm (#1)}]}
\def\qtn{\begin{equation}}
\def\qtnl#1{\begin{equation}\label{#1}}
\def\eqtn{\end{equation}}
\def\lmm{\begin{lemma}}
\def\lmml#1{\begin{lemma}\label{#1}}
\def\elmm{\end{lemma}}
\def\crllr{\begin{corollary}}
\def\crllrl#1{\begin{corollary}\label{#1}}
\def\ecrllr{\end{corollary}}
\def\hpthss{\begin{hypothesis}}
\def\hpthssl#1{\begin{hypothesis}\label{#1}}
\def\ehpthss{\end{hypotxesis}}
\def\xmpl{\begin{example}}
\def\xmpll#1{\begin{example}\label{#1}}
\def\exmpl{\end{example}}
\def\css{\begin{cases}}
\def\ecss{\end{cases}}

\begin{document}
\title[Coset closure of a circulant S-ring and schurity problem]{Coset closure of a circulant S-ring\\ and schurity problem} 
\author{Sergei Evdokimov}
\address{St. Petersburg Department of Steklov Institute of Mathematics, Russia}
\email{evdokim@pdmi.ras.ru}
\author{Ilya Ponomarenko}
\address{St. Petersburg Department of Steklov Institute of Mathematics, Russia}
\thanks{The second author was partially supported by RFFI Grant 14-01-00156}
\email{inp@pdmi.ras.ru}
\maketitle

\begin{abstract}
Let $G$ be a finite group. There is a natural Galois correspondence between the permutation groups containing $G$
as a regular subgroup, and the Schur rings (S-rings) over~$G$. The problem we deal with in the paper, is to 
characterize those S-rings that are closed under this correspondence, when the group $G$ is cyclic
(the schurity problem for circulant S-rings). It is proved that up to a natural reduction,
the characteristic property of such an S-ring  is to be a certain algebraic fusion of its coset closure introduced
and studied in the paper. Basing on this characterization we show that the schurity problem is equivalent 
to the consistency of a modular linear  system associated with a circulant S-ring under consideration. As a byproduct  we 
show that a circulant S-ring is  Galois closed if and only if so is its dual.
\end{abstract}

\section{Introduction}
A {\it Schur ring} or {\it S-ring} over a finite group $G$ can be defined as a subring of the group ring $\mZ G$
that is a free $\mZ$-module spanned by a partition of $G$ closed under taking inverse and containing $\{1_G\}$
as a class (see Section~\ref{140214a} for details). It is well known that there is
a Galois correspondence~\footnote{We recall that a Galois correspondence between two posets
consists of two mappings reversing the orders such that both superpositions are closure operators.}
between the permutation groups on $G$ that contain the regular group $G_{right}$, and the S-rings over~$G$:
\qtnl{030314a}
\{\Gamma\le\sym(G):\ \Gamma\ge G_{right}\}\ \rightleftarrows\ 
\{\cA\le\mZ G:\ \cA\ \text{is an S-ring over}\ G\}.
\eqtn
More precisely, the "$\rightarrow$"~mapping is given by taking the partition of $G$ into the orbits 
of the stablizer  of $1_G$ in $\Gamma$, whereas the "$\leftarrow$"~mapping is 
given by taking the automorphism group of the colored Cayley graph corresponding to the partition of~$G$ 
associated with~$\cA$.
The Galois closed objects are called $2$-closed 
groups and {\it schurian} S-rings, respectively. The schurity problem consists in finding an inner characterization of
schurian S-rings.\medskip

The theory of S-rings was initiated by I.~Schur (1933) and later developed by H.~Wielandt~\cite{Wie64} and his followers.
The starting point for Schur was  the Burnside theorem stating that any primitive permutation 
group containing a regular cyclic $p$-group of composite order, is 2-transitive. Using the S-ring method
introduced by him, Schur generalized this theorem to an arbitrary finite cyclic group  $G$ (cf.~Theorem~\ref{040214a}). 
To some extent this explains the fact that "Schur had conjectured for a long time that every S-ring over~$G$ is 
determined by a suitable permutation group"~\cite[p.54]{Wi69}.
This statement  had been known as the Schur-Klin conjecture up to 2001, when the first examples
of circulant (i.e.~over a cyclic group) S-rings  were constructed in~\cite{EP01} by the authors. A recent
result in \cite{EP11} shows that schurian circulant S-rings are relatively rare. In this paper
we provide a solution to the schurity problem for circulant S-rings.\medskip

The non-schurian examples of S-rings were constructed using the operation of {\it generalized wreath product} 
introduced in~\cite{EP01} (and independently in~\cite{LM96} under the name "wedge product").
This  is not suprising due to the seminal Leung-Man theorem according to which any circulant S-ring 
can be constructed from S-rings of rank~$2$ and cyclotomic S-rings by means of two operations: tensor product and 
generalized wreath product \cite{LM96}. Here under a {\it cyclotomic} S-ring $\cA$
we mean the ring of all $K$-invariant elements of $\mZ G$ where $K$ is a subgroup of $\aut(G)$: 
\qtnl{190214d}
\cA=(\mZ G)^K.
\eqtn
The Leung-Man theorrm reduces the schurity problem  for circulant S-rings to finding a criterion for the schurity of the generalized wreath product. Such a criterion, based on a generalization of the Leung-Man theory (see~\cite{EP01ce}),
was obtained in paper~\cite{EP12} where the generalized wreath product of permutation groups was introduced
and studied. All these results form a background to prove the main results of the paper 
(see Sections~\ref{140214a}--\ref{280114a}).\medskip

Let $\cA$ be a circulant S-ring. Suppose that among the "bricks" in the Leung-Man decomposition of~$\cA$,
there is a non-cyclotomic S-ring. Then this ring is  of rank~$2$, its underlying group has composite
order and it is Cayley isomorphic to the restriction of~$\cA$  to one of its sections. Moreover, 
as it was proved in~\cite{EP12} the S-ring~$\cA$
has a quite a rigid structure that enables us to control the schurity of~$\cA$. This provides a reduction of the
schurity problem to the case when $\cA$ has no rank~$2$ section of composite order. The S-rings
satisfying the latter property are {\it quasidence} in sense of paper~\cite{EP11} (Theorem~\ref{060214a}). 
Thus without loss of generality we concentrate on the schurity problem for quasidence S-rings.\medskip

Our first step is to represent the schurian closure $\sch(\cA)$
of a quasidense circulant S-ring $\cA$ in a regular form (Theorem~\ref{251212a}).
The idea here is to replace the ring~$\cA$ by a simpler one keeping the structure of its Leung-Man decomposition.
The simplification is achieved by changing each "brick" for a group ring. This leads to the class of {\it coset} 
S-rings, i.e. ones for  which any class of the corresponding partition of the group $G$ 
is a coset of a subgroup in~$G$.
It appears that this class is closed under restriction to a section, tensor and generalized wreath products, 
and consists of schurian quasidense S-rings (Theorems~\ref{170611b} and~\ref{150113a}).  
The regular form of $\sch(\cA)$we want to come, will be defined by means of the following concept.

\dfntnl{201113a}
The coset closure of a quasidense circulant S-ring $\cA$ is the intersection $\cA_0$ of all coset S-rings over~$G$ 
that contain~$\cA$. 
\edfntn

The coset closure of any quasidense circulant S-ring is a coset S-ring (Theorem~\ref{311013a}).
Now, to  clarify how to represent the schuran closure of $\cA$ via its coset closure, suppose that
the group $G$ is of prime order. In this case it is well known that the S-ring $\cA$  is of the 
form~\eqref{190214d}, and, moreover, $\cA_0=\mZ G$. In particular, $\cA$ is  schurian and
any automorphism of $G$ induces a similarity of $\cA_0$.\footnote{Under a similarity of an S-ring $\cA$ we mean a ring 
isomorphism of it  that respects the partition of~$G$ corresponding to $\cA$, see Subsection~\ref{190214e}.}
Furthermore,  if the automorphism belongs to the group $K$, the similarity is  identical on~$\cA$.
Thus
$$
\cA=(\mZ G) ^K=(\cA_0)^{\Phi_0}
$$
where $\Phi_0=\Phi_0(\cA)$ is the group of all similarities of $\cA_0$ that are identical on~$\cA$. It appears that
this idea works for any quasidense S-ring~$\cA$.

\thrml{251212a}
Let $\cA$ be a quasidense circulant S-ring. Then
$$
\sch(\cA)=(\cA_0)^{\Phi_0}.
$$
In particular, $\cA$ is schurian if and only if $\cA=(\cA_0)^{\Phi_0}$.
\ethrm

Theorem~\ref{251212a} gives a necessary and sufficient condition for an S-ring to be schurian. This condition
being a satisfactory from the theoretical point of view, is hardly an inner characterization. To obtain the
latter, we prove Theorem~\ref{170113a} below. Let us discuss briefly the idea behind it.\medskip

One of the key properties of coset S-rings that is used in the proof of Theorem~\ref{251212a}, is that
every similarity of any such ring is induced by isomorphism (Theorem~\ref{150113a}).
This fact also shows that in the schurian case the set of all isomorphisms of $\cA_0$ that induce similarities belonging 
to~$\Phi_0$, 
forms a permutation group the associated S-ring  of which coincides with~$\cA$. In general, this  is not true.
A rough reason for this can be explained as follows. Set
\qtnl{110214b}
\frS_0=\{S\in\frS(\cA_0):\ (\cA_0)_S=\mZ S\}
\eqtn
where $\frS(\cA_0)$ is the set of all $\cA_0$-sections and $(\cA_0)_S$ is the restriction of $\cA_0$ to~$S$. Then in 
the schurian case every S-ring $\cA_S$ with $S\in\frS_0$,   must be cyclotomic,
whereas in general
this condition does not necessarily hold. However, if even all the S-rings $\cA_S$ are cyclotomic,
one still might find a section~$S$ for which $\cA_S\ne\sch(\cA)_S$.  These two reasons are controlled respectively by conditions~(1) and~(2) of Theorem~\ref{170113a}.\medskip

It should be mentioned that the proof  of the fact that the circulant S-rings constructed in~\cite{EP01} are non-schurian,
was based on studying the relationship between their cyclotomic sections. More careful analysis can be found 
in~\cite{EP03} where   the isomorphism problem for circulant graphs was solved. In that paper the authors introduce
and study the notion of projective  equivalence on the sections of a circulant S-ring (this notion is similar to one used in 
the lattice theory). It appears that the class $\frS_0$ defined in~\eqref{110214b} is closed under the projective 
equivalence and taking subsections (Corollary~\ref{221113a}).  Moreover,
$$
\frS(\cA_0)=\frS(\cA)
$$ 
(statement~(1) of Theorem~\ref{281212b}). \medskip 

To formulate Theorem~\ref{170113a} we need additional notation. For  $S\in\frS(\cA)$ 
denote by $\aut_\cA(S)$
the subgroup of  $\aut(S)$ that consists of all Cayley automorphisms of the S-ring~$\cA_S$.
A family 
$$
\fS=\{\sigma_S\}_{S\in\frS_0}
$$ 
is called a {\it multiplier} of~$\cA$ if for any sections $S,T\in\frS_0$
such that $T$ is projectively equivalent to a subsection of~$S$ the automorphisms $\sigma_T\in\aut(T)$ and
$\sigma_S\in\aut(S)$ are induced by raising to the same power\footnote{We recall that any
automorphism of a finite cyclic group is induced by raising to a power coprime to the order of this group.}.
The set of all multipliers of $\cA$ forms a subgroup  of the direct product $\prod_{S\in\frS_0}\aut_\cA(S)$
that is denoted by $\mult(\cA)$.

\thrml{170113a}
A quasidense circulant S-ring $\cA$ is schurian if and only if  the following two conditions are satisfied for all $S\in\frS_0$:
\nmrt
\tm{1} the S-ring $\cA_S$ is cyclotomic,
\tm{2} the restriction homomorphism from $\mult(\cA)$ to $\aut_\cA(S)$ is surjective.
\enmrt
\ethrm

By Theorem~\ref{291112u}  the class $\frS_0$ consists of all $\cA$-sections~$S$ such
that each Sylow subgroup of $S$ (treated as a section of~$G$) is projectively equivalent to a subsection 
of a principal $\cA$-section. Thus in contrast to Theorem~\ref{251212a}, Theorem~\ref{170113a} gives
a necessary and sufficient condition for an S-ring~$\cA$ to be schurian in terms of $\cA$ itself rather than  of 
its coset closure~$\cA_0$. It should be remarked that in general the  class $\frS_0$  may contain non-cyclotomic
sections. However, we do not know whether condition~(1) in Theorem~\ref{170113a} is implied by
condition~(2).\medskip

In fact, the starting point of this paper was  the following question: 
is the property of an S-ring "to be schurian" preserved under taking the dual. 
The following theorem  deduced in Section~\ref{100214b} from Theorem~\ref{170113a},
answers this question in the positive. 

\thrml{170113b}
A circulant S-ring is schurian if and only if so is the S-ring dual to it.
\ethrm

We would like to reformulate  Theorem~\ref{170113a} in the number theoretical language. 
In what follows we assume that condition~(1)  of that theorem is satisfied. To make condition~(2) more clear 
let us fix a section $S_0\in\frS_0$ and an integer $b$ coprime to $n_{S_0}=|S_0|$ for which the mapping 
$s\mapsto s^b$, $s\in S_0$, belongs to $\aut_\cA(S_0)$. Let us consider
the following system of linear equations in integer variables $x_S$, $S\in\frS_0$:
\qtnl{091013a}
\css
x_S\equiv x_T\pmod{n_T} ,              &\\
x_{S_0}\equiv b\pmod{n_{S_0}}   &\\
\ecss
\eqtn
where $S$ and $T$ run over $\frS_0$  and the section $T$ is projectively equivalent to a subsection 
of~$S$. We are interested only in the solutions of this system that satisfy the additional condition
\qtnl{120214a}
(x_S,n_S)=1\quad\text{for all}\quad S\in\frS_0.
\eqtn
Every such solution produces the family $\fS=\{\sigma_S\}$ 
where $\sigma_S$ is the automorphism of the group $S$ taking~$s$ to~$s^{x_S}$. Moreover, the equations 
in the first line of  \eqref{091013a} guarantee that if a section $T$ is projectively equivalent to a subsection 
of~$S$, then the automorphisms $\sigma_T\in\aut(T)$ and
$\sigma_S\in\aut(S)$ are induced by raising to the same power. Therefore, 
$$
\fS\in\mult(\cA).
$$
Conversely, it is easily seen that given $S_0\in\frS$ every multiplier of~$\cA$ 
produces a solution of system~\eqref{091013a} for the corresponding~$b$. Finally, the  consistency
of this system for all~$S_0$ and all possible~$b$ is equivalent to  the surjectivity of the restriction homomorphism from $\mult(\cA)$ to $\aut_\cA(S_0)$ for
all~$S_0$. Thus we come to  the following corollary of Theorem~\ref{170113a}.

\crllrl{091013b}
Let $\cA$ be a quasidense circulant S-ring  such that for any section  $S\in\frS_0$,  the S-ring $\cA_S$ is cyclotomic.
Then $\cA$ is schurian if and only if system~\eqref{091013a} has a solution satisfying~\eqref{120214a}
for all possible $S_0$ and $b$.
\ecrllr

Corollary~\ref{091013b} reduces the schurity problem for circulant S-rings to solving modular 
linear system~\eqref{091013a}
under restriction~\eqref{120214a}. One possible way to solve this system is to represent the group
$\prod_S\aut_\cA(S)$ as a permutation group on the disjoint union of the sections~$S$. Then every equation
in the first line of~\eqref{091013a} defines a subgroup of that group the index of which is at most $n^2$.
Therefore the set of solutions can be found by a standard permutation group technique, see \cite[p.~144]{L93} for details.\medskip

Concerning permutation groups we refer to~\cite{DM}. For the reader convenience an extended background
on the S-ring theory including new material,   is given in Sections~\ref{140214a}--\ref{280114a}. In 
Sections~\ref{241013a} and~\ref{200214a} we study liftings of
generalized wreath products in the non-dense and dense cases respectively. The theory of coset S-rings
is developed in Sections~\ref{051113a} and~\ref{200214b}, and culminates in Theorem~\ref{150113a} showing
that these S-rings are schurian and separable. The coset closure and multipliers of a quasidense S-ring are introduced and
studied in Sections~\ref{111113b} and~\ref{200214c} respectively. Theorems~\ref{251212a} and~\ref{170113a}
are proved in Section~\ref{200214f}. In the final Section~\ref{100214b} we prove Theorem~\ref{170113b}.\medskip

{\bf Notation.}
As usual by $\mZ$ we denote the ring of rational integers.

The set of all right cosets of a subgroup $H$ in a group $G$   is denoted by $G/H$.
For a set $X\subset G$ we put $X/H=\{Y\in G/H:\ Y\subset X\}$.

The subgroup of $G$ generated by a set $X\subset G$ is denoted by $\lg X\rg$; we also set
$
\rad(X)=\{g\in G:\ gX=Xg=X\}
$.

For a prime $p$ a Sylow $p$-subgroup of~$G$ is denoted by $G_p$.


For a section $S=U/L$ of $G$ the quotient epimorphism from $U$ onto $S$ is denoted by $\pi_S$.

For a bijection $f:G\to G'$ and a set $X\subset 2^G$ or an element $X\in 2^G$, the induced bijection from $X$ onto 
$X^f$ is denoted by~$f^X$.

The group of all permutations of $G$ is denoted by $\sym(G)$.

For a set $\Delta\subset\sym(G)$ and a section $S$ of  $G$ we set
$$
\Delta^S=\{f^S:\ f\in \Delta,\ S^f=S\}.
$$

The subgroup of $\sym(G)$ induced by right multiplications of ~$G$  is denoted by $G_{right}$.

For $X,Y\in G/H$ we set $G_{X\to Y}=\{f^X:\ f\in G_{right},\ X^f=Y\}$.

The holomorph $\hol(G)$ is identified with the subgroup of~$\sym(G)$ generated by $G_{right}$ and~$\aut(G)$.

The orbit set of a group $\Gamma\le\sym(G)$ is denoted by $\orb(\Gamma)=\orb(\Gamma,G)$.

We write $\Gamma\twoe\Gamma'$ if groups $\Gamma,\Gamma'\le\sym(G)$ are $2$-equivalent, i.e.
have the same orbits in the coordinate-wise action on~$G\times G$.

For a positive integer $n$, the cyclic group the elements of which are integers modulo~$n$, is denoted 
by  $\mZ_n$.

For an automorphism $\sigma$ of a cyclic group $G$ of order $n$ the unique
element $m\in\mZ_n$ for which $x^\sigma=x^m$, $x\in G$, is denoted by~$m(\sigma)$.

\section{S-rings}\label{140214a}
In what follows we use the  notation and terminology from paper~\cite{EP12} where a part of the material
is contained.  Concerning the basic S-ring theory and duality theory we also refer to~\cite[Ch.2.6]{BI}
and~\cite[Ch.4]{Wie64}, respectively.\medskip

\sbsnt{Definitions and basic facts.}
Let $G$ be a finite group with identity $1_G$. A subring~$\cA$ of the group ring~$\mZ G$ is called a 
{\it Schur ring} ({\it S-ring}, for short) 
over~$G$ if  there exists a partition $\cS(\cA)$ of~$G$ such that
\nmrt
\tm{1} $\{1_G\}\in\cS(\cA)$,
\tm{2} $X\in\cS(\cA)\ \Longrightarrow\ X^{-1}\in\cS(\cA)$,
\tm{3} $\cA=\Span\{\sum_{x\in X}x:\ X\in\cS(\cA)\}$.
\enmrt
Two S-rings $\cA$ and $\cA'$ are called {\it Cayley isomorphic} if there exists a ring isomorphism from $\cA$ onto $\cA'$
that is induced by isomorphism of underlying groups; the latter  is called {\it Cayley isomorphism} from $\cA$ 
onto $\cA'$. Obviously, $\cS(\cA)^f=\cS(\cA')$ for any such isomorphism~$f$.  \medskip

The elements of $\cS(\cA)$ and the number $\rk(\cA)=|\cS(\cA)|$ are called respectively the {\it basic sets} and the {\it rank} of the 
S-ring~$\cA$. Any union of basic sets is called an {\it $\cA$-subset of~$G$} or {\it $\cA$-set}. The set of all of them is 
closed with respect to taking inverse and product, and forms a lattice with respect to inclusion. Given an $\cA$-set $X$ the 
submodule of~$\cA$ spanned by the set
$$
\cS(\cA)_X=\{Y\in\cS(\cA):\ Y\subset X\}
$$
is denoted by $\cA_X$. \medskip

A subgroup of $G$ that is an $\cA$-set is called  {\it $\cA$-subgroup} of~$G$ or {\it $\cA$-group}; the set of all of them 
is denoted by $\cG(\cA)$. With each $\cA$-set $X$ one can  naturally associate two $\cA$-groups, namely $\lg X\rg$ and  
$\rad(X)$ (see Notation). The S-ring $\cA$ is called {\it dense} if every subgroup of~$G$ is an $\cA$-group, and {\it primitive} 
if the only $\cA$-subgroups are $1$ and $G$.\medskip

A  section $S=U/L$ of the group $G$ is  called a section of~$\cA$ or an  {\it $\cA$-section}, if both  $U$ and $L$ are 
$\cA$-groups. In this case we also say that  $\cA$ contains~$S$. The module
$$
\cA_S=\Span \{\pi_S(X):\ X\in\cS(\cA)_U\}
$$
is an S-ring over the group~$S$, the basic sets of which are exactly the sets in the right-hand side of the formula.
A section $S$ is of rank~$2$ (resp. primitive, cyclotomic) if so is the S-ring $\cA_S$. The set of all (resp. all cyclotomic) 
$\cA$-sections is denoted by $\frS(\cA)$ (resp. $\frS_{cyc}(\cA)$). Set
$$
\cF_{prin}(\cA)=\{\grp{X}/\rad(X):\ X\in\cS(\cA)\}.
$$
Any element of this set is called a {\it principal} $\cA$-section.\medskip

The partial order on the set of all S-rings over $G$ that is induced by inclusion, is denoted by~$``\le"$. Thus 
$\cA_1\le\cA_2$ if and only if any basic  set of $\cA_1$ is a union of basic  sets of~$\cA_2$
(in this case we say that $\cA_2$ is an extension of~$\cA_1$). 
The least and greatest elements are 
$$
\Span\{1,\sum_{x\in G}x\}\qaq\mZ G
$$ 
respectively. Next, the module $\cA_1\cap\cA_2$  is an S-ring, the basic sets of which 
form the finest partition of $G$ that is coarser than both $\cS(\cA_1)$ and $\cS(\cA_2)$. 
It follows that the  set of all S-rings over $G$ forms a lattice in which the meet (resp. join) of
$\cA_1$ and~$\cA_2$ coincides with $\cA_1\cap\cA_2$ (resp.  with the intersection of all S-rings over $G$ that 
contain both~$\cA_1$ and~$\cA_2$). It is easily seen that 
\qtnl{311013b}
(\cA_1\cap\cA_2)_S=(\cA_1)_S\cap(\cA_2)_S
\eqtn
for any section $S\in\frS(\cA_1)\cap\frS(\cA_2)$.\medskip

Let $K\le\aut(G)$ be a group of Cayley isomorphisms of the S-ring~$\cA$. Then the set $\cA^K$ of the
elements in $\cA$ left fixed under the induced action of $K$ on~$\mZ G$, forms an S-ring over~$G$. Any such
ring with $\cA=\mZ G$, is called {\it cyclotomic}  and is denoted by $\cyc(K,G)$. The classes of the
corresponding partition of $G$ are exactly the orbits of the group~$K$.\medskip

If $\cA_1$ and $\cA_2$ are S-rings over groups $G_1$ and $G_2$ respectively, then the subring $\cA=\cA_1\otimes \cA_2$ 
of the ring $\mZ G_1\otimes\mZ G_2=\mZ G$ where $G=G_1\times G_2$, is an S-ring over the group $G$ with
$$
\cS(\cA)=\{X_1\times X_2: X_1\in\cS(\cA_1),\ X_2\in\cS(\cA_2)\}.
$$
It is called the {\it tensor product} of $\cA_1$ and $\cA_2$.

\sbsnt{Generalized wreath product.}
Let $S=U/L$ be a section of an S-ring~$\cA$.  We say that~$\cA$ is an {\it $S$-wreath product} if the group $L$ is 
normal in~$G$ and $L\le\rad(X)$ for all basic sets $X$  outside~$U$; in this case we write
\qtnl{050813a}
\cA=\cA_U\wr_S\cA_{G/L},
\eqtn
and omit $S$ when $|S|=1$; in the latter case $\cA$ is called {\it wreath product}. When an explicit 
reference to the section~$S$ is not important, we use the term {\it generalized wreath 
product}. The $S$-wreath product is {\it nontrivial} or {\it proper}  if $1\ne L$ and $U\ne G$.

\thrml{230114a}
Let $S=U/L$ be a section of a group~$G$, and let  $\cA_1$ and $\cA_2$ be S-rings over 
the groups $U$ and $G/L$ respectively such that $S$ is both an $\cA_1$- and an $\cA_2$-section with
$$
(\cA_1)_S=(\cA_2)_S.
$$
Then the set of S-rings $\cA$ such that $\cA_U=\cA_1$ and $\cA_{G/L}=\cA_2$ has the smallest element. Moreover,
it is a unique $S$-wreath product in this set.
\ethrm
\proof It was proved in \cite[Theorem~3.1]{EP01}  that the $S$-wreath product $\cA$ from~\eqref{050813a} is uniquely determined
and belongs to the set of S-rings from the theorem statement. For any other S-ring $\cA'$ from this set the preimage in $\cS(\cA')$
of at least one basic set in $\cS(\cA_2)$ outside $U/L$ contains at least two basic sets, whereas the same preimage in $\cS(\cA)$
consists of one element. This proves that $\cA'\ge\cA$, and hence the minimality of~$\cA$.\bull\medskip

\sbsnt{Duality.}\label{091213a}

Let $\cA$ be an S-ring over a finite abelian group~$G$ and $\wh G$ the
group dual to~$G$, i.e. the group of all irreducible complex characters of~$G$.
Given $S\subset G$ and $\chi\in\wh G$ set
\qtnl{201108b}
\chi(S)=\sum_{s\in S}\chi(s).
\eqtn
Characters $\chi_1,\chi_2\in\wh G$ are called equivalent if
$\chi_1(S)=\chi_2(S)$ for all $S\in\cS(\cA)$. Denote by $\wh\cS$ the set of
classes of this equivalence relation.
Then the submodule of $\mZ \wh G$ spanned by the elements $\xi(X)$, $X\in\wh\cS$,
is an S-ring over~$\wh G$ (see~\cite[Theorem~6.3]{BI}). This ring is called
{\it dual} to~$\cA$ and is denoted by~$\wh\cA$. Obviously, $\cS(\wh\cA)=\wh\cS$.
Moreover, $\rk(\wh\cA)=\rk(\cA)$ and
\qtnl{311008a}
\cG(\wh\cA)=\{H^\bot:\ H\in \cG(\cA)\}
\eqtn
where $H^\bot=\{\chi\in\wh G:\ H\le\ker(\chi)\}$. It is also true that the S-ring
dual to $\wh\cA$ is equal to~$\cA$.\medskip

It is easily seen that the mapping from $\aut(G)$ to $\aut(\wh G)$ that takes $\sigma$ to $\wh\sigma$ defined by
$\chi^{\wh\sigma}(g)=\chi(g^\sigma)$, is a group isomorphism. The image of a group $K\le\aut(G)$ is denoted by 
$\wh K$. 

\thrml{210114a}
Let $\cA=\cyc(K,G)$ where $K\le \aut(G)$. Then $\wh\cA=\cyc(\wh K,\wh G)$.
\ethrm
\proof Let $X\in\orb(\wh K,\wh G)$. Then given $\chi_1,\chi_2\in X$ there
exists $\sigma\in K$ such that $\chi^{}_1=\chi_2^{\,\wh\sigma}$. Since $S=S^\sigma$ for
each basic set $S$ of~$\cA$, this implies that
$$
\chi^{}_1(S)=\chi_2^{\,\wh\sigma}(S)=\chi^{}_2(S^\sigma)=\chi^{}_2(S),\qquad S\in\cS(\cA).
$$
Therefore 
$
\cyc(\wh K,\wh G)\ge\wh\cA  
$. Replacing here $K$ and $G$  by $\wh K$ and $\wh G$ respectively, we conclude that the S-ring $\cA$
contains the S-ring dual to $\cyc(\wh K,\wh G)$. By duality this implies that
$
\wh\cA\ge\cyc(\wh K,\wh G)$.
Thus $\wh\cA=\cyc(\wh K,\wh G)$.\bull\medskip 

Some more properties of the dual S-ring are contained in the following statement
proved in~\cite[Theorems~2.4, 2.5]{EP13}. In what follows given a section $S=U/L$ the group
$\wh S$ is canonically identified with the section $L^\bot/U^\bot$ of the group~$\wh G$ that is
called the section {\it dual} to~$S$. In particular, if $G=G_1\times G_2$, then
$\wh G=\wh{G_1}\times\wh{G_2}$.


\thrml{140509a}
Let $\cA$ be an S-ring over an abelian group $G$. Then
\nmrt
\tm{1} $\wh{\cA_S}=\wh\cA_{\wh S}$ for any $S\in\frS(\cA)$,
\tm{2} $\cA=\cA_1\otimes\cA_2$ if and only if $\wh\cA=\wh\cA_1\otimes\wh\cA_2$,
\tm{3} $\cA$ is an $S$-wreath product if and only if $\wh\cA$ is an $\wh S$-wreath product.\bull
\enmrt
\ethrm

\section{Similarities, isomorphisms and schurity}

In this section we follow papers \cite{EP01,EP01ce,EP12} except for terminology: similarities and isomorphisms defined
below were called in~\cite{EP01,EP01ce} weak and strong isomorphisms respectively
(see also the remark in Subsection~\ref{170114r}). The reader familiar with association scheme theory will see
that the definitions of similarity, isomorphism, etc. given in this section, are compatible with those used
for Cayley schemes (see~\cite{EP12}).

\sbsnt{Similarities.}\label{190214e}
Let $\cA$ and $\cA'$ be S-rings over groups $G$ and $G'$ respectively. A ring isomorphism $\varphi:\cA\to\cA'$ is called
{\it similarity} from $\cA$ to $\cA'$, if for any $X\in\cS(\cA)$ there exists $X'\in\cS(\cA')$ such that
$$
\varphi(\sum_{x\in X}x)=\sum_{x'\in  X'}x'.
$$
It follows from the definition that the mapping $X\mapsto X'$
is a bijection   from~$\cS(\cA)$ onto~$\cS(\cA')$. This bijection is naturally extended to a 
bijection between $\cA$- and $\cA'$-sets, that takes $\cG(\cA)$ to $\cG(\cA')$, and hence 
$\frS(\cA)$ to $\frS(\cA')$. The images of an $\cA$-set~$X$ and $\cA$-section~$S$ are denoted by $X^\varphi$ 
and $S^\varphi$ respectively. For any such $S$ the similarity $\varphi$ induces a similarity
$$
\varphi_S:\cA^{}_{S^{}}\to\cA'_{S'}\quad\text{where}\quad S'=S^\varphi.
$$
 The set of all similarities from $\cA$ to $\cA'$ is denoted by $\Phi(\cA,\cA')$;
we also set $\Phi(\cA)=\Phi(\cA,\cA)$.\medskip

The above bijection between the $\cA$- and $\cA'$-sets is in fact an isomorphism of the corresponding lattices. 
It follows that given an $\cA$-set $X$ we have
\qtnl{150114e}
\grp{X^\varphi}=\grp{X}^\varphi\qaq\rad(X^\varphi)=\rad(X)^\varphi.
\eqtn
These equalities together with the obvious equation $X=X\rad(X)$, immediately imply the following statement.

\lmml{221113f}
Any similarity of an S-ring is uniquely determined by its restrictions to principal sections.\bull
\elmm

Any automorphism $\sigma\in\aut(G)$ can be extended linearly
to a ring automorphism $\varphi_\sigma$ of $\mZ G$. Thus the lemma below immediately follows
from the definition of similarity.

\lmml{271213a}
The mapping $\sigma\mapsto\varphi_\sigma$ is a group isomorphism from $\aut(G)$ onto $\Phi(\mZ G)$.\bull
\elmm

Let $\Phi$ be a group of similarities of the S-ring~$\cA$. Then the set $\cA^\Phi$ of the
elements in $\cA$ left fixed under the action of $\Phi$, is obviously an S-ring over~$G$ for which
$$
\cS(\cA^\Phi)=\{X^\Phi:\ X\in\cS(\cA)\}
$$
where $X^\Phi=\bigcup_{\varphi\in\Phi}X^\varphi$. When $\cA=\mZ G$, from Lemma~\ref{271213a} it follows that 
$\cA^\Phi$ is a cyclotomic S-ring.\medskip

We complete the subsection by describing the similarities of generalized wreath products; mostly
this was done in~\cite{EP01}.

\thrml{271213x}
Let $\cA$ and $\cA'$ be S-rings over abelian groups $G$ and $G'$. Suppose that 
$\cA$ is an $S$-wreath product where $S=U/L$. Then 
\nmrt
\tm{1} for any similarity $\varphi\in\Phi(\cA,\cA')$ the S-ring $\cA'$ is the $S'$-wreath product where 
$S'=S^\varphi$,
\tm{2} if $\cA'$ is an $S'$-wreath product,  then the mapping $\varphi\mapsto(\varphi_U,\varphi_{G/L})$
induces a bijection from the set $\{\varphi\in\Phi(\cA,\cA'):\ S'=S^\varphi\}$  to the set
$$
\{(\varphi_1,\varphi_2)\in\Phi(\cA^{}_{U^{}},\cA'_{U'})\times \Phi(\cA^{}_{G^{}/L^{}},\cA'_{G'/L'}):
(\varphi_1)_S=(\varphi_2)_S\}.
$$
where $U'=U^\varphi$ and $L'=L^\varphi$.
\enmrt
\ethrm
\proof Statement~(2) follows from \cite[Theorem~3.3]{EP01} . 
To prove statement~(1) it suffices to verify that $L'X'=X'$ for all $X'\in\cS(\cA')_{G'\setminus U'}$. However, 
this is true because $(G\setminus U)^\varphi=G'\setminus U'$ and $LX=X$ for all $X\in\cS(\cA)_{G\setminus U}$.\bull

\sbsnt{Isomorphisms.}\label{170114r}
Let $\cA$ and $\cA'$ be S-rings over groups $G$ and $G'$ respectively. A bijection $f:G\to G'$ is called
an {\it isomorphism} from $\cA$ onto~$\cA'$ if there exists a similarity $\varphi\in\Phi(\cA,\cA')$ such that
given $X\in\cS(\cA)$ we have
\qtnl{260613a}
f(Xy)=X^\varphi f(y)\quad \text{for all}\quad y\in G,
\eqtn
or, equivalently, $f(x)f(y)^{-1}\in X^\varphi$ for all $x,y\in G$ with $xy^{-1}\in X$. In this case we also say
that $f$ {\it induces}~$\varphi$. Clearly, any isomorphism induces a uniquely determined similarity. 
The set of all isomorphisms and isomorphisms with a fixed~$\varphi$
are denoted by $\iso(\cA,\cA')$ and $\iso(\cA,\cA',\varphi)$ respectively.\medskip

It follows from the definition that the isomorphism $f$ that takes $1_{G^{}}$ to $1_{G'}$, takes
$\cS(\cA)$ to $\cS(\cA')$ and satisfies the condition $f(Xy)=f(X)f(y)$ for all $X\in\cS(\cA)$ and $y\in G$. Therefore
$f$ is a strong isomorphism from $\cA$ to $\cA'$ in the sense of~\cite{EP01}. Conversely, according to that paper
any strong isomorphism $f$ induces a similarity $\varphi$ such that $f(X)=X^\varphi$, and hence satisfies~\eqref{260613a}.
Thus S-rings are isomorphic if and only if they are strongly isomorphic (see equality~\eqref{150114a} below).\medskip

It immediately follows from the definition that
\qtnl{150114a}
G^{}_{right}\iso(\cA,\cA',\varphi)G'_{right}=\iso(\cA,\cA',\varphi).
\eqtn
In particular, not every isomorphism takes $1_{G^{}}$ to $1_{G'}$. But even if it does, it is not necessarily a Cayley isomorphism.
However, $|G|=|G'|$ because every similarity preserves the order of the underlying group. Moreover, 
since equality~\eqref{260613a} obviously holds also for any $\cA$-set $X$, and, in particular, 
every isomorphism preserves the right cosets of any $\cA$-group.

\lmml{220513a}
In the above notation  let $H$ be an $\cA$-group and $H'=H^\varphi$. Then
$$
hf^X h'\in\iso(\cA^{}_{H^{}},\cA'_{H'},\varphi^{}_H)
$$
for all $X\in G/H$, $h\in G^{}_{H^{}\to X^{}}$ and $h'\in G'_{X'\to H'}$ where  $X'=X^f$.
\elmm
\proof Denote by $g$ and $g'$  the permutations from $G^{}_{right}$ and $G'_{right}$ such that $g^H=h$ and $(g')^{X'}=h'$, 
respectively. Then by~\eqref{150114a} the bijection $gfg':G\to G'$ induces the similarity $\varphi$. So the bijection
$hf^X h'=(gfg')^H$ induces the similarity $\varphi_H$ as required.\bull\medskip

The following statement characterizes the isomorphisms of a generalized wreath product.

\thrml{140513a}
Let $\cA$ and $\cA'$ be S-rings over abelian groups $G$ and $G'$ and $\varphi:\cA\to\cA'$ a similarity. Suppose 
that $\cA$ is a $U/L$-wreath product. Then the set  $\iso(\cA,\cA',\varphi)$ consists of all bijections $f:G\to G'$ such 
that $(G/U)^f=G'/U'$, $(G/L)^f=G'/L'$ where $U'=U^\varphi$ and $L'=L^\varphi$,  and
\qtnl{250613b}
f^{G/L}\in \iso(\cA^{}_{G/L},\cA'_{G/L},\varphi^{}_{G/L}),
\quad
gf^Xg'\in \iso(\cA^{}_U,\cA'_U,\varphi^{}_U)
\eqtn
for all $X\in G/U$ and some $g\in G_{U\to X}$ and $g'\in G'_{X'\to U'}$ where $X'=X^f$.
\ethrm
\proof Set $F=\iso(\cA,\cA',\varphi)$ and denote by $F'$ the set of all bijections $f:G\to G'$ 
satisfying~\eqref{250613b}. Then the inclusion $F\supset F'$
immediately  follows from the basic properties of similarities and Lemma~\ref{220513a}. Conversely, let  $f\in F'$.
We have to verify that equality~\eqref{260613a} holds for all $X\in\cS(\cA)$ and $y\in G$. Suppose first that $X\subset G\setminus U$.
Then from the equality $(G/L)^f=G'/L'$ and the first relation in~\eqref{250613b} it follows that
$$
f^{G/L}(X^\pi y^\pi)=(X^\pi)^{\varphi_{G/L}} f^{G/L}(y^\pi)
$$
where $\pi=\pi_{G/L}$ and $y\in G$. After taking the preimages of both sides in the latter equality
we obtain that $f(XyL)=(XL)^\varphi f(yL)$. On the other hand, $L\le\rad(X)$. Due to~\eqref{150114e} this implies that
$L'\le\rad(X^\varphi)$. Thus
$$
f(Xy)=f(XyL)=(XL)^\varphi f(yL)=X^\varphi f(y)L'=X^\varphi f(y),
$$
which proves the required statement in our case.
Let now $X\subset U$ and $y\in G$. Then from the equality $(G/U)^f=G'/U'$ and the second relation in~\eqref{250613b} it follows that
$$
X^{gf^{Y}g'}=X^\varphi
$$
for some $g\in G^{}_{U^{}\to Y^{}}$ and $g'\in G'_{Y'\to U'}$ where $Y=Uy$ and $Y'=Y^f$. By Lemma~\ref{220513a} without loss 
of generality we can assume that the permutations $g$ and $g'$ are induced by multiplications by $y$ and $(y^f)^{-1}$
respectively.
Then $X^{gf^Yg'}=(Xy)^f(y^f)^{-1}$, and hence
$$
f(Xy)=(Xy)^f=X^\varphi y^f=X^\varphi f(y)
$$
and we are done.\bull\medskip

Let $\cK$ be a class of S-rings closed under Cayley isomorphisms.
Following~\cite{EP09a}  an S-ring $\cA$ is called separable with respect to $\cK$ if
$\iso(\cA,\cA',\varphi)\ne\emptyset$ for all similarities $\varphi:\cA\to\cA'$ where
$\cA'\in\cK$. In this paper we say that a circulant S-ring is {\it separable} if it is
separable with respect to the class of all circulant S-rings.

\sbsnt{Schurity.}

Let $G$ be a finite group. It was proved by Schur (see \cite[Theorem~24.1]{Wie64}) that  any group
$\Gamma\le\sym(G)$ that contains $G_{right}$ produces an S-ring~$\cA$ over $G$ such that
$$
\cS(\cA)=\orb(\Gamma_1,G)
$$
where $\Gamma_1=\{\gamma\in\Gamma:\ 1^\gamma=1\}$ is the stabilizer of the point $1=1_G$ in $\Gamma$. Any such  
S-ring is called {\it schurian}.  Group rings and S-rings of rank~$2$ are obviously schurian. \medskip

Let $\Gamma$ and $\Delta$ be permutation groups on $G$ such that $G_{right}\le\Gamma\cap\Delta$. Then it is easily
seen that the S-ring associated with the group $\grp{\Gamma,\Delta}$ equals the intersection of S-rings associated
with $\Gamma$ and $\Delta$. It follows that the intersection of schurian S-rings is schurian. Therefore so is the S-ring
\qtnl{240214a}
\sch(\cA)=\bigcap_{\cA'\ge\cA,\\ \cA'\ \text{is schurian}}\cA'.
\eqtn
It is called the {\it schurian closure} of~$\cA$. Clearly, $\sch(\cA)\ge\cA$, and the equality is attained  if and only if 
the S-ring $\cA$ is schurian.\medskip

The schurity concept is closely related to  automorphisms of an S-ring. In contrast to a common algebraic tradition
the automorphism group of an S-ring $\cA$ is not defined to be $\iso(\cA,\cA)$; in accordance with a combinatorial
tradition we set $\aut(\cA)=\iso(\cA,\cA,\id_\cA)$. Thus $f\in\aut(\cA)$  if and only if 
\qtnl{291113a}
f(Xy)=Xf(y),\qquad X\in\cS(\cA),\quad y\in G. 
\eqtn
The latter is equivalent to say that given $X\in\cS(\cA)$ we have $f(x)f(y)^{-1}\in X$ whenever $xy^{-1}\in X$.
Therefore $\aut(\cA)$ can also be defined as the automorphism group of the colored Cayley graph corresponding
to the partition $\cS(\cA)$ of the group~$G$ (cf. Introduction).\medskip

It follows that any basic set of $\cA$ is invariant with respect to the group $\aut(\cA)_1$, whereas 
equality~\eqref{150114a} shows that $G_{right}\le\aut(\cA)$.
Moreover, the group $\aut(\cA)$ is the largest subgroup of $\sym(G)$ that satisfies these two properties.  
Now, let $\Gamma$ be a permutation group the schurian S-ring $\sch(\cA)$ is associated with. Then since 
$\sch(\cA)\ge\cA$,
the maximality  of $\aut(\cA)$ implies that 
$$
\Gamma\le\aut(\cA).
$$ 
It follows that $\sch(\cA)$ contains the S-ring associated with $\aut(\cA)$. Therefore 
these S-rings are equal (see~\eqref{240214a}). In fact, this shows that the closure of the S-ring~$\cA$
with respect to Galois correspondence~\eqref{030314a} coincides with $\sch(\cA)$. Thus
the above definition of a schurian S-ring  is compatible with that given in Introduction.\medskip


Let $f\in\aut(\cA)_1$.
Then any $\cA$-set (in particular, $\cA$-group)  is invariant with respect to the automorphism~$f$. 
Moreover,
for any $\cA$-section $S$ we have $f^S\in\aut(\cA_S)$. In particular, the S-ring $\cA_S$ is schurian whenever so is~$\cA$.\medskip 

The following result proved in \cite[Corollary~5.7]{EP12}  gives a criterion for the schurity of generalized wreath products 
that will be  repeatedly used throughout the paper. Below we set
$$
\cM(\cA)=\{\Gamma\le\sym(G):\ \Gamma\twoe\aut(\cA)\ \,\text{and}\ \,G_{right}\le\Gamma\}.
$$

\thrml{060813a}
Let $\cA$ be an S-ring over an abelian group~$G$. Suppose that $\cA$ is an $S$-wreath product where~$S=U/L$. Then $\cA$
 is schurian if and only if so are the S-rings $\cA_{G/L}$ and $\cA_U$ and there exist groups $\Delta_0\in\cM(\cA_{G/L})$ and 
 $\Delta_1\in\cM(\cA_U)$ such that $(\Delta_0)^S=(\Delta_1)^S$.
\ethrm


It should be remarked that a permutation $f\in\sym(G)$ preserving every basic set of $\cA$ does not necessarily 
belong to $\aut(\cA)$. However, as the lemma below shows, this is so when, for example, $f\in\aut(G)$.

\lmml{021213a}
Let $\cA$ be an S-ring and $\cA'\ge\cA$. Then the group of all isomorphisms
of $\cA'$ that fix every element of~$\cA$, is a subgroup of $\aut(\cA)$.
\elmm
\proof Let $f$ be an isomorphism of $\cA'$ and $\varphi$ the similarity of $\cA'$ induced by~$f$.
Then $f(X'y)=(X')^\varphi f(y)$ for all $\cA'$-sets $X'$ and $y\in G$. Since $\cA'\ge\cA$, this is true 
for all $X'\in\cS(\cA)$. So if $f$ fixes every element of~$\cA$, then $\varphi$ is identical on~$\cA$, and hence
$f(Xy)=Xf(y)$ for all $X\in\cS(\cA)$. Thus $f\in\aut(\cA)$.\bull


\section{Sections in S-rings}

\sbsnt{Projective equivalence.}\label{300114a}
Let $G$ be a group. Denote by $\cG(G)$ the set of its subgroups, and by $\frS(G)$ the set of its sections, 
i.e. quotients of subgroups of~$G$.
When this does not lead to misunderstanding, we write $H$ instead of $H/1$ where $H\in\cG(G)$, and identify $\frS(S)$ with
the corresponding subset of $\frS(G)$ where $S\in\frS(G)$.\medskip

A section $U'/L'$ is called a {\it subsection} of a section $U/L$ if $U'\le U$ and $L'\ge L$; in this case we write
$U'/L'\preceq U/L$. This defines a partial order on the set  $\frS(G)$. This order has the greatest element $G/1$; 
any minimal element is of the form $H/H$ where $H\in\cG(G)$.\medskip

A section $U/L$ is called a {\it multiple} of a section $U'/L'$ if
\qtnl{310114b}
L'=U'\cap L\quad\text{and}\quad U=U'L.
\eqtn
The {\it projective equivalence} relation "$\sim$"  on the set $\frS(G)$ is defined to be the transitive closure of the relation "to be a multiple". 
Any two projecively equivalent sections are obviously isomorphic as groups.
The set of all equivalence classes is denoted  by $\cP(G)$. Under a {\it quasisubsection} of a section~$S$ we mean any section which
is projectively equivalent to a subsection of~$S$.\medskip

Given a section $S=H/K$ of a group $G$ we define a surjection $\rho_{G,S}$ from the subgroups of~$G$ to
the subgroups of~$S$  by
$$
\rho_{G,S}(U)=(U\cap H)K/K.
$$
Let us  extend this mapping to $\cF(G)$ by $\rho_{G,S}(U/L)=\rho_{G,S}(U)/\rho_{G,S}(L)$. Using the above
identification of $\frS(S)$ with the corresponding subset of $\frS(G)$, we obtain that
\qtnl{170613c}
\rho_{G,S}(U/L)=(U\cap H)K/(L\cap H)K.
\eqtn
Clearly, this mapping is identical on the set $\cF(S)$. Moreover, $\rho_{G,S}$ inducess a surjective homomorphism
of the corresponding partially ordered sets. To simplify notations we will write $T_S$ instead of $\rho_{G,S}(T)$
when the group $G$ is fixed.\medskip

If the group $G$ is abelian, then the set $\cG(G)$ of all subgroups of $G$ forms a modular lattice
\footnote{We recall that a lattice is modular, if $x\vee(y\wedge z)=(x\vee y)\wedge z$ whenever $x\le z$.} 
in which the join  and meet of $H$ and $K$ are defined as $HK$ and $H\cap K$ respectively~\cite{S94}.

\lmml{210613a}
Let $G$ be an abelian group. Then given $S,T\in\cF(G)$ we have $S_T\sim T_S$.
\elmm
\proof Let $S=H/K$ and $T=U/L$. Then by the definition of $S_T$ and $T_S$ we have
$$
S_T=(H\cap U)L/(K\cap U)L\qaq  T_S=(U\cap H)K/(L\cap H)K.
$$
A straightforward check shows that both $S_T$ and $T_S$ are multiples of the section $(H\cap U)/(H\cap U\cap KL)$ 
(we made use the fact the lattice $\cG(G)$ is modular). Therefore  the sections $S_T$ and $T_S$ are projectively 
equivalent.\bull\medskip



\sbsnt{Restrictions.}
Let $\cA$ be an S-ring over a group~$G$. Set
$$
\cP(\cA)=\{C\cap\frS(\cA):\ C\in\cP(G),\ C\cap\frS(\cA)\ne\emptyset\}.
$$
Then $\cP(\cA)$ forms a partition of the set $\frS(\cA)$ into classes  of projectively
equivalent $\cA$-sections. It should be noted that if $\cA'\ge\cA$, then $\frS(\cA')\supset\frS(\cA)$ and each class of projectively equivalent $\cA$-sections is contained
in a unique class of projectively equivalent $\cA'$-sections.

\thrml{261010a}
Let $\cA$ be an S-ring over a group $G$. Then given projectively equivalent $\cA$-sections $S$ and $T$
there exists a Cayley isomorphism~$f$ from $\cA_S$ onto $\cA_T$ such that $(\gamma^S)^f=\gamma^T$ for all 
$\gamma\in\aut(\cA)$ leaving the point~$1_G$ fixed.
\ethrm
\proof Follows from \cite[Lemma~3.1,Theorem~3.2]{EP12}.\bull\medskip

Obviously, any two sections of a class $C\in\cP(G)$ have the same order; we call it the {\it order} of this class. If, in addition, 
$C\in\cP(\cA)$ where $\cA$ is an S-ring over~$G$, then from Theorem~\ref{261010a} it follows that all sections in~$C$ have 
the same rank~$r$, and also if $C$ contains a primitive (resp. cyclotomic, dense)  section, then all sections 
in~$C$ are primitive (resp. cyclotomic, dense). In these cases we say that~$C$ is a class of rank~$r$, and a primitive 
(resp. cyclotomic, dense) class.\medskip

In general, an $\cA$-section projectively equivalent to a principal $\cA$-section is not principal.
However, at least in the cyclic group case such a section is
{\it subprincipal}, i.e. a subsection of a principal $\cA$-section (Lemma~\ref{100214a}).
But even in this case the class of subprincipal sections is not closed under the projective equivalence.
To avoid this inconvenience we define an 
$\cA$-section to be {\it quasiprinciple} (resp. {\it quasisubprinciple}) if it is projectively equivalent to  a principal 
(resp. subprincipal) $\cA$-section.\medskip

Let $S\in\frS(\cA)$. Then the mapping $\rho_{G,S}:T\mapsto T_S$ defined in Subsection~\ref{300114a}
induces a mapping from $\cF(\cA)$ to $\cF(\cA_S)$ that is denoted 
by the same letter. The following statement shows that it preserves generalized wreath products.

\thrml{031212c}
Let $\cA$ be an S-ring over an abelian group $G$. Suppose that $\cA$ is the $T$-wreath product where $T\in\frS(\cA)$. Then 
$\cA_S$ is the $T_S$-wreath product for all~$S\in\frS(\cA)$.
\ethrm
\proof Let $T=U/L$ and $S=H/K$. We have to prove that $L_S\le\rad(Y)$ for all $Y\in \cA_{S\setminus U_S}$.
However, any $Y\in\cS(\cA_S)$ is of the form $XK/K$ for some $X\in\cS(\cA_H)$. If, in addition, $Y\in \cA_{S\setminus U_S}$,
then $X\in \cA_{G\setminus U}$, and hence $L\le \rad(X)$ because $\cA$ is the $T$-wreath product.
Since $\rad(X)\le\grp{X}\le H$, this implies that 
$$
L_S\le\rad(X)_S=\rad(X)K/K\le\rad(XK/K)=\rad(Y)
$$ 
as required.\bull

\sbsnt{Duality.} Let $G$ be an abelian group.  Then the mapping $H\mapsto H^\bot$ induces a lattice
antiisomorphism from $\cG(G)$ onto $\cG(\wh G)$ (see e.g.~\cite{S94}).

\lmml{310114a}
For any sections $S,T\in\frS(G)$ the following statements hold:
\nmrt
\tm{1} $S\sim T$  if and only if $\wh S\sim\wh T$; moreover,
$T$ is a multiple of $S$ if and only if $\wh S$ is a multiple of $\wh T$,
\tm{2} $\wh{\rho_{G,S}(T)}=\rho_{\wh G,\wh S}(\wh T)$,
\tm{3} $S\preceq T$ if and only if $\wh S\preceq\wh T$.
\enmrt
\elmm
\proof Let $S=H/K$ and $T=U/L$. To prove statement~(1) without loss of generality we can assume that
$T$ is a multiple of $S$. Then from~\eqref{310114b} it follows that
$$
K^\bot=H^\bot L^\bot\quad\text{and}\quad U^\bot=H^\bot\cap L^\bot,
$$
which means that $\wh S=K^\bot/H^\bot$ is a multiple of $\wh T=L^\bot/U^\bot$ as required. Statement~(2) 
follows from~\eqref{170613c} and the modularity of the lattice $\cG(G)$:
$$
\wh{\rho_{G,S}(T)}=((U\cap H)K)^\bot/((L\cap H)K)^\bot=((UK)\cap H)^\bot/((LK)\cap H)^\bot=
$$
$$
(U^\bot\cap K^\bot) H^\bot/(L^\bot\cap K^\bot)H^\bot=\rho_{\wh G,\wh S}(\wh T).
$$
Statement~(3) is obvious.\bull\medskip

Let $\cA$ be an S-ring over $G$. Given a class $C\in\cP(\cA)$ we define the dual class by
$$
\wh C=\{\wh S:\ S\in C\}.
$$
Then from statement (1) of Lemma~\ref{310114a} it follows that $\wh C\in \cP(\wh\cA)$. Moreover,
the classes $\wh C$ and $C$ have the same order and rank, and if one of them is primitive (resp. cyclotomic, dense),
then so is the other one (see equality~\eqref{311008a} and Theorem~\ref{210114a}).

\section{Circulant S-rings}\label{280114a}

\sbsnt{General theory.}
We begin with a well-known result on circulant primitive S-rings, that goes back to  Burnside and Schur. Despite
the fact that they dealt with groups, their results can be interpreted as results on schurian circulant
S-rings. Moreover, the Schur method works in the non-schurian case as well
(see, e.g.~\cite{EP05}).

\thrml{040214a}
Any circulant primitive S-ring is of rank~$2$, or a cyclotomic S-ring over a group of prime order.
\ethrm

Let $\cA$ be an S-ring over a cyclic group~$G$. Then 
$$
X^\sigma\in\cS(\cA)\quad\text{for all}\quad X\in\cS(\cA),\ \sigma\in\aut(G),
$$
see \cite[Theorem~23.9]{Wie64}. Let now $X$ be a {\it highest} basic
set of~$\cA$, i.e. one containing a generator of~$G$. Then the above statement implies that the group $\rad(X)$
does not depend on the choice of~$X$. It is called  the {\it radical} of~$\cA$ and denoted
by $\rad(\cA)$. The following statement proved by Leung and Man is a cornerstone of the 
circulant S-ring theory (see \cite[Corollaries~5.5,6.4]{EP01ce}).

\thrml{141113g}
Let $\cA$  be a circulant S-ring. Then the following two statements hold:
\nmrt
\tm{1} $\rad(\cA)\ne 1$ if and only if $\cA$ is a proper generalized wreath product,
\tm{2} $\rad(\cA)= 1$ if and only if $\cA$ is the tensor product of a cyclotomic S-ring with trivial
radical and S-rings of rank~$2$.\bull
\enmrt
\ethrm

Any principal $\cA$-section $S$ is obviously a section with trivial radical in the sense that $\rad(\cA_S)=1$.
The converse is not true, but we have the following statement.

\lmml{100214a}
In a circulant S-ring, any section with trivial radical is subprincipal.
\elmm
\proof Suppose that $\rad(\cA_S)=1$ where $S=U/L$ is a section of a circulant S-ring $\cA$.  
We observe that every highest basic set $X$ of the S-ring~$\cA_U$ produces a highest basic
set $\pi(X)$ of the S-ring $\cA_S$ where $\pi=\pi_S$. Thus
$$
\pi(\rad(X))\le\rad(\pi(X))=\rad(\cA_S)=1.
$$
It follows that $\rad(X)\le L$ . Therefore $S\preceq T$ where $T=U/\rad(X)$. Since $U=\grp{X}$,
the $\cA$-section $T$ is principal, and we are done.\bull\medskip

From Theorem~\ref{261010a} it follows that the set of all $\cA$-sections with trivial radical
is closed with respect to the projective equivalence. Another property of the projective equivalence for
circulant S-rings is given in the following statement proved in~\cite[Lemma~5.2]{EP03}.
Below  a section of a class $C\in\cP(\cA)$  is called the {\it smallest} (respectively, {\it largest}) one if every section 
of $C$, is a multiple of it (respectively, if it is a multiple of every section of~$C$). 

\thrml{040214b}
Any class of projectively equivalent sections of a circulant S-ring
has the largest and smallest elements.\bull
\ethrm

\sbsnt{Quasidence S-rings.}
A circulant S-ring $\cA$ is called {\it quasidense}, if any primitive
$\cA$-section is of prime order. Any dense S-ring is obviously quasidense.  Moreover, in the quasidense case
any minimal $\cA$-group is of prime
order, any maximal $\cA$-group is of prime index, and the S-ring $\cA_S$ is
dense for any $\cA$-section $S$ of prime power order.

\thrml{060214a}
A circulant S-ring is quasidense if and only if it contains no rank~$2$ sections of composite order.
\ethrm
\proof By Theorem~\ref{040214a} a section of  composite order is of rank~$2$  if and only if
it is primitive. Thus the required statement immediately follows from the
definition of a quasidence S-ring.\bull\medskip

It is easily seen that the class of quasidense S-rings is closed
under restriction to a section, and the tensor and generalized wreath products.

\thrml{180214a}
Any extsension of quasidense circulant S-ring is quasidense.
\ethrm
\proof Let $\cA'$ be a non-quasidense extension of a quasidense S-ring~$\cA$ over a cyclic group $G$. Then by Theorem~\ref{060214a} the S-ring $\cA'$ contains a rank~$2$ section $S=U/L$ of composite order. 
Denote by $H$ an $\cA$-group  of prime order (such  a group 
does exist because $\cA$ is quasidense). We claim that
\qtnl{190214a}
H\cap U=1\quad \text{or}\quad H\le L.
\eqtn 
Indeed, if~\eqref{190214a} is not true, then $H\le U$ and $H\cap L=1$ because the
order of $H$ is prime. It follows that the S-ring $\cA'_{S}$ contains the group $HL/L$ of prime order $|H|$. 
The latter group does not equal $S$ because the order of $S$ is composite. However this is impossible because 
$\rk(\cA'_S)=2$.\medskip

By~\eqref{190214a} the S-ring $\cA'_{G/H}$ contains the section $S'=UH/LH$ which is projectively equivalent to the section~$U/L$ by the modularity of the lattice~$\cG(G)$. By Theorem~\ref{261010a} this implies that
$$
\rk(\cA'_{S'})=\rk(\cA'_S)=2.
$$
Thus, $\cA'_{G/H}$ is a non-quasidense extension of the quasidense 
S-ring $\cA_{G/H}$ due to Theorem~\ref{060214a}. Assuming without loss of generality that the order of the group $G$ is minimal possible we come to a contradiction.\bull

\thrml{060214b}
The class of quasidense circulant S-rings with trivial radical is closed with respect to taking extensions, and
consists of cyclotomic, and hence dense S-rings.
\ethrm
\proof Let $\cA$ be a quasidense circulant S-ring with trivial radical. Then it is cyclotomic, and hence dense
by \cite[Theorem~3.1]{EP11}. If now $\cA'\ge\cA$, then the S-ring $\cA'$ is dense, and hence quasidense. 
Suppose  on the contrary that
$\rad(\cA')\ne 1$. Then by statement~(1) of Theorem~\ref{141113g} the S-ring $\cA'$ is a proper 
$S'$-wreath product for some $\cA'$-section~$S'$. This section is also $\cA$-section by the density
of~$\cA$. Therefore $\cA$ is the $S'$-wreath product, which is impossible by statement~(1) of 
Theorem~\ref{141113g} because $\rad(\cA)=1$.\bull\medskip

The following theorem deduced from~\cite[Theorem~3.5]{EP11} shows that any schurian quasidense circulant S-ring
can be obtained from an appropriate solvable permutation group that "locally"  has a rather 
simple form. 

\thrml{230312a}
Let $\cA$ be a schurian quasidense circulant S-ring. Then there exists a group $\Gamma\in\cM(\cA)$ such that
$\Gamma^S=\hol_\cA(S)$ for any quasisubprinciple $\cA$-section $S$ 
where $\hol_\cA(S)=\hol(S)\cap\aut(\cA_S)$.
\ethrm
\proof By \cite[Theorem~3.5]{EP11} there exists a group $\Gamma\in\cM(\cA)$ such that
$\Gamma^S=\hol_\cA(S)$ for any $\cA$-section $S$ such that $\rad(\cA_S)=1$. Since
obviously $\hol_\cA(S)^T=\hol_\cA(T)$ for all $\cA$-sections $T\preceq S$, the required
statement follows from Theorem~\ref{261010a} for $\gamma\in\Gamma$.\bull

\sbsnt{Duality.}
The following two theorems establishing selfdual properties of a circulant S-ring will be used
in Section~\ref{100214b} to prove Theorem~\ref{170113b}.

\thrml{201213a}
Let $\cA$ be a circulant S-ring. Then
\nmrt
\tm{1} $\rad(\cA)=1$ if and only if $\rad(\wh\cA)=1$,
\tm{2} $\cA$  is quasidense  if and only if $\wh\cA$ is quasidence.
\enmrt
\ethrm
\proof Statement~(1) immediately follows from statement~(3) of Theorem~\ref{140509a} and
Theorem~\ref{141113g}. Statement~(2) holds because the primitivity  and the order of an $\cA$-section
are preserved under duality: the former by equality~\eqref{311008a} and statement~(1) of Theorem~\ref{140509a},
whereas the latter by the definition of the dual section.\bull

\thrml{101213e}
A section of a circulant S-ring is  quasisubprincipal  if and only if so is its dual. 
\ethrm
\proof Let $S$ be a quasisubprincipal section of a circulant S-ring~$\cA$. Then $S$ is projectively equivalent
to a subsection $S'$ of a principal $\cA$-section $T$. Without loss of generality we can assume that $T$
is a $\preceq$-maximal principal $\cA$-section. Then  $\wh T$ is a $\preceq$-maximal $\wh\cA$-section
by  statement~(1) of Lemma~\ref{310114a}. Moreover,
by statement~(2) of Theorem~\ref{201213a} we have $\rad(\wh\cA_{\wh T})=1$. Thus by
Lemma~\ref{100214a} we conclude that
$\wh T$ is a principal $\wh\cA$-section. On the other hand, from statements~(1) and~(3) of Lemma~\ref{310114a}
we have $\wh S\sim \wh{S'}$ and  $\wh{S'}\preceq{\wh T}$. Thus $\wh S$ is a quasisubprincipal $\cA$-section
as required. The converse statement holds by duality.\bull

\section{Lifting: nondense case}\label{241013a}

We begin with a characterization of the dense circulant S-rings in terms of forbidden subsections. To do this
we will say that an S-ring $\cA$ over a cyclic group $G$ is {\it elementary nondense} if $|G|$ is a composite number and
$\cA$  has rank~$2$, or $|G|$ is the product of two distinct primes
and $\cA$ is a proper wreath product. In the former case the S-ring is not quasidense, whereas in the latter case it is.  A section $S$ of a 
circulant S-ring $\cA$ is called {\it elementary nondense} if the S-ring $\cA_S$ is elementary nondense.  

\lmml{151013a}
The section that is dual or projectively equivalent to an elementary nondense section is elementary nondense.
\elmm
\proof Follows from Theorems~\ref{140509a}  and~\ref{261010a}.\bull\medskip

It is easily seen that a circulant S-ring that contains an elementary nondense  section cannot be dense. This proves 
the ``if" part of the following statement.

\thrml{301012a}
A circulant S-ring~$\cA$ is not dense if and only if there exists an elementary nondense $\cA$-section.
In particular, any minimal nondense $\cA$-section is elementary.
\ethrm
\proof To prove the ``only if" part suppose that $\cA$ is not dense. If it is not quasidense, then 
it contains a rank~$2$ section $S$ of composite order (Theorem~\ref{060214a}). Since
$S$ is elementary nondense, we are done. Suppose that $\cA$ is quasidesnse and $S=V/K$ is a 
minimal nondense $\cA$-section.  Then there exists  a non-$\cA$-group $H\le G$ such that
$$
K<H<V.
$$
The minimality of $S$ implies that $K$ is a maximal $\cA$-group inside~$H$. So after decreasing $H$ (if necessary)
without loss of generality we can assume that the number $p=|H/K|$ is prime. Next, by the quasidensity of $\cA_{V/K}$   
there exists an $\cA$-group~$M$ such that $K\le M\le V$ and  the number
$q=|M/K|$ is prime. Since $H$ is not an $\cA$-group, it follows that $p\ne q$ and hence
$$
H\cap M=K.
$$
Moreover, by the minimality of $S$, the section $V/M$ is dense. Therefore $MH$ is an $\cA$-group. So $MH=V$ by
the minimality of $S$. Thus  $|V/K|=pq$.  To complete the proof it suffices to
note that a quasidence but not dense S-ring over a cyclic group of order $pq$ has to be a proper wreath 
product (Theorem~\ref{141113g}).\bull\medskip

In what follows $\cA$ is a quasidense  S-ring over a cyclic group $G$.\footnote{Generalized wreath products arising in nonquasidense case
had been studied in~\cite{EP12}.} Suppose that $S_0=V/K$ is an elementary nondense $\cA$-section. 
Denote by $H$ the unique $\cA$-subgroup of $G$ such that $K<H<V$. Then
$$
\cA_{S_0}=\cA_{H/K}\wr\cA_{V/H}.
$$
The largest $\cA$-section which is projectively equivalent to $V/H$, is obviously of the form
$UV/U$ for some $\cA$-group $U=U(S_0)$; similarly, the smallest $\cA$-section which is projectively equivalent to $H/K$, is of the form $L/(K\cap L)$ for some $\cA$-group $L=L(S_0)$. (The existence of the largest and smallest
sections follows from Theorem~\ref{040214b}.) Clearly,
\qtnl{290113c}
1<L\le U<G.
\eqtn
The relevant part of the lattice of $\cA$-groups is given in Fig.~\ref{f8}.\medskip 

The definition of the section $U/L$ associated with $S_0$ is uniform in the following sense.
Let $S$ be an $\cA$-section  that contains $S_0$ as a subsection. Since the mapping $\rho_{G,S}$
defined in Subsection~\ref{300114a} induces a lattice epimorphism from $\cG(\cA)$ to~$\cG(\cA_S)$, we have
\qtnl{191113a}
\rho_{G,S}(U/L)=U_S/L_S
\eqtn
where $U_S/L_S$ is the section of the S-ring $\cA_S$ defined by $S_0$ in the same way as the section
$U/L$ in the S-ring~$\cA$. In particular, $U_{S_0}/L_{S_0}=H/H$.\medskip

\def\VRT#1{*=<8mm>[o][F-]{#1}}
\begin{figure}[t]
\tiny{
$\hspace{45mm}\xymatrix@R=10pt@C=20pt@M=0pt@L=5pt{
& \VRT{UV}\ar@{-}[dr]\ar@{-}[dl] & & & \\
\VRT{U} \ar@{-}[dr] & & \VRT{V}\ar@{-}_q[dl] & &\\
& \VRT{H}\ar@{-}[dr]\ar@{-}_{p}[dl] & &  & \\
\VRT{K}\ar@{-}[dr] & & \VRT{L}\ar@{-}[dl] & &  \\
& \VRT{K\cap L} & &  & \\
}$
}
\caption{}\label{f8}
\end{figure}

The following two statements will be used to prove the main result of this section (Theorem~\ref{290113a}).

\lmml{310113a}
In the above notation the following statements hold:
\nmrt
\tm{1} the section~$S_0$  is a quasisubsection  of any  $\cA$-group $M\not\le U$,
\tm{2} the section~$S_0$ is a quasisubsection  of  any $\cA$-section $G/N$ with $N\not\ge L$.
\enmrt
\elmm
\proof From Lemma~\ref{151013a} it follows that  $\wh{S_0}$ is an elementary nondense section of the S-ring $\wh\cA$. Therefore statement~(2)
follows from statement~(1) by duality. To prove statement~(1) we claim that $UV/U$ is a $\cA$-quasisubsection of~$M$.
 Indeed, since $U\le UM\cap UV\le UV$ and the number $|UV/U|=|V/H|$ is prime, we have
$$
UM\cap UV\in\{U,UV\}.
$$
However, $UM\cap UV\ne U$, because otherwise $UVM/UM$ is obviously a multiple of the section $UV/U$, which contradicts the maximality of it.
Therefore $UM\cap UV=UV$, and hence $UV\le UM$. But $U(M\cap UV)=UM\cap UV$ because the lattice $\cG(\cA)$
is modular. Then
$$
U(M\cap UV)=UM\cap UV=UV.
$$
Since also obviously $U\cap M\cap UV=U\cap M$, we conclude that the section $UV/U$ is a multiple of $(M\cap UV)/(M\cap U)$.
This proves the claim because the latter is an $\cA_M$-section.\medskip

 Set $V'/H'$ to be the smallest $\cA$-section which is projectively equivalent to $V/H\sim UV/U$. Then by the claim in the
previous paragraph, $V'/H'$ is an $\cA_M$-section. To complete the proof it suffices to verify that
$$
V'/K'\sim V/K
$$
where $K'=K\cap V'$. Let us prove that $V/K$ is a multiple of $V'/K'$. Suppose on the contrary that $KV'\ne V$. 
Then $KV'\le H$ because $K\le KV'\le V$ and $H$ is the only $\cA$-group strictly between $K$ and $V$.
Therefore $V'\le H$. On the other hand, $V=V'H$ because $V/H$ is a multiple of $V'/H'$. Thus
 $V=V'H\le H$ which is impossible.\bull

\crllrl{290113b}
The section $S_0$ is a quasisubsection of any $\cA$-section $M/N$ with $M\not\le U$ and $L\not\le N$.
\ecrllr
\proof By statement~(1) of Lemma~\ref{310113a}  the S-ring $\cA_M$ contains a section $S'_0=V'/K'$ projectively equivalent to~$S_0=V/K$.
Since the section $S'_0$ is elementary nondense (Lemma~\ref{151013a}), there is a unique $\cA$-group  $H'$ strictly between
$K'$ and $V'$. From the choice of  the group $L$ it follows that
$$
H'/K'\sim H/K\sim L/K\cap L,
$$
Moreover, $L\le H'\le M$. Since also $N\le M$, the hypothesis of Lemma~\ref{310113a} is satisfied for
$G=M$ and $S^{}_0=S'_0$. So by statement~(2) of this lemma we conclude that the section $M/N$ has 
a subsection projectively equivalent to $S'_0\sim S^{}_0$ as required.\bull\medskip

From the Leung-Man theory it follows that any quasidense circulant S-ring that is not dense, is a proper generalized wreath product
(see~\cite[Theorem~5.3]{EP01ce}). The following theorem gives an explicit form of such a product.

\thrml{290113a}
Let $\cA$ be a quasidense circulant S-ring and $S_0$ an elementary nondense section. Then $\cA$ is a proper $U/L$-wreath product
with $U=U(S_0)$ and $L=L(S_0)$. Moreover, any dense $\cA$-section is either an $\cA_U$-section or $\cA_{G/L}$-section.
\ethrm
\proof To prove the first statement  it suffices to verify by~\eqref{290113c} that 
$$
\rad(X)\ge L\quad\text{for all}\quad X\in\cS(\cA)_{G\setminus U}.
$$
Suppose on the contrary that this is not true for some~$X$. Then the hypothesis of Corollary~\ref{290113b} holds for
$M=\grp{X}$ and $N=\rad(X)$. Therefore $S_0$ is a quasisubsection of~$M/N$.
Since $S_0$ is elementary nondense, this implies that the S-ring $\cA_{M/N}$ is not dense. However, this is impossible because
$\cA_{M/N}$ is a quasidense S-ring with trivial radical (Theorem~\ref{060214b}).\medskip

 To prove the second statement suppose on the contrary that there exists
a dense $\cA$-section $M/N$ which is neither an $\cA_U$- nor $\cA_{G/L}$-section. Then $S_0$ is a quasisubsection of $M/N$
by Corollary~\ref{290113b}. Therefore the S-ring $\cA_{M/N}$ is not dense. Contradiction.\bull\medskip

The following auxiliary statement will be used in Section~\ref{111113b}.

\lmml{111113c}
Let $\cA'$ be an extension of a quasidense circulant S-ring~$\cA$. Then $\cH(\cA')=\cH(\cA)$
if and only if any elementary nondense section of~$\cA$ is an elementary nondense section of~$\cA'$. 
\elmm
\proof The ``only if'' part is obvious. To prove the ``if'' part suppose on the contrary that there exists a group $H\in\cH(\cA')\setminus\cH(\cA)$. Without loss of generality we 
assume that the cyclic group $G$ underlying $\cA$ and $\cA'$ is minimal possible. Then $G\ne 1$. Moreover,
\qtnl{111113e}
U\in\cH(\cA)\ \ \&\ \ U\ne G\quad \Longrightarrow\quad U\cap H=1
\eqtn
and
\qtnl{111113f}
L\in\cH(\cA)\ \ \&\ \ L\ne 1\quad \Longrightarrow\quad LH=G.
\eqtn
Indeed, relation~\eqref{111113e} follows from relation~\eqref{111113f}  by duality. To prove~\eqref{111113f}
suppose on the contrary that $LH\ne G$. Since $L$ and $H$ are $\cA'$-groups, so is the group~$LH$. By the
minimality of $G$ we have $\cH(\cA^{}_{G/L})=\cH(\cA'_{G/L})$, and hence $LH$ is an $\cA$-group.  Again by
the minimality this implies that $\cH(\cA'_{LH})=\cH(\cA^{}_{LH})$. It follows that $H$ is an $\cA$-group. 
Contradiction.\medskip

By the quasidensity of $\cA$ there exist $\cA$-groups $U$ and $L$ such that the numbers $|G/U|$ and
$|L|$ are prime. Since $1<H<G$, from~\eqref{111113e} and \eqref{111113f} it follows that $G=U\times H=L\times H$.
Thus the numbers $|H|$ and $|G/H|$ are prime. Therefore the S-ring $\cA$ is elementary nondense.
By the lemma hypothesis this implies that so is $\cA'$. This implies that  $H\not\in\cH(\cA')$ in
contrast to the choice of~$H$.\bull

\section{Lifting: dense case}\label{200214a}

In this section we are to get an analog of the theory developed in Section~\ref{241013a} but this time for dense S-rings. In what follows
under a {\it $p$-section} we mean a section which is a $p$-group.

\thrml{241013b}
Let $\cA$ be a circulant dense S-ring. Then $\rad(\cA)\ne 1$ if and only if $\cA$ contains  a non-quasisubprinciple $p$-section.
\ethrm
\proof The "if" part is obvious. To prove the "only if" part suppose that $\rad(\cA)\ne 1$. Then by 
\cite[Theorem~5.4]{EP01ce} the S-ring~$\cA$ is a $U/L$-wreath product such that 
\qtnl{110314a}
|G/U|=|L|=p
\eqtn
where $p$ is a prime and $G$ is the underlying group of $\cA$. To complete the proof
it suffices to verify that the $\cA$-section~$G_p$ is non-quasisubprinciple. Suppose on the contrary that $G_p$ is a quasisubsection of a
principal $\cA$-section $T=H/K$. Then $|T_p|=|G_p|$ and hence
\qtnl{251013c}
K_p=(G/H)_p=1.
\eqtn
On the other hand, since $\cA$ is the $U/L$-wreath product, Theorem~\ref{031212c} implies that the S-ring $\cA_T$ is a 
$U_T/L_T$-wreath product where 
$$
U_T/L_T=\rho_{G,T}(U/L)=(U\cap H)K/(L\cap H)K.
$$
Using~\eqref{110314a} and~\eqref{251013c} we obtain by comparing $p$-parts  that $(U\cap H)K<H$ 
and  $(L\cap H)K>K$. It follows that
$U_T\ne T$ and $L_T\ne 1$. Therefore $\cA_T$ is a proper generalized wreath product. However, this is impossible
because the section $T$ is principle and hence $\rad(\cA_T)=1$.\bull\medskip

As the following example shows, not every minimal non-quasi\-subprin\-ciple section is a $p$-section.\medskip

{\bf Example.}\label{230414a}
 Let $p,q,r$ be distinct primes such that $r-1$ is divided by both $p$ and $q$. Let us define 
S-rings $\cA_1$ and $\cA_2$ over the group $\mZ_{rp^2q^2}$ as follows
$$
\cA_1=\cyc(K_{1,1},\mZ_{rp^2})\otimes \cyc(K_{1,2},\mZ_{q^2}),
$$
$$
\cA_2= \cyc(K_{2,1},\mZ_{p^2})\otimes\cyc(K_{2,2},\mZ_{rq^2})
$$
where $K_{1,1}$, $K_{1,2}$, $K_{2,1}$ and $K_{2,2}$  are groups of order $p$, $q$, $p$ and $q$ respectively. The groups $K_{1,2}$ and $K_{2,1}$
are uniquely determined; choose the groups $K_{1,1}$ and $K_{2,2}$ so that none of the coordinate projections is trivial. Then
$\cyc(K_{1,1},\mZ_{rp^2})$ and $\cyc(K_{2,1},\mZ_{rq^2})$ are S-rings with trivial radicals. Moreover,
using natural identifications we have
$$
(\cA_1)^{p^2q^2}=\cyc(K_{2,1},\mZ_{p^2})\otimes  \cyc(K_{1,2},\mZ_{q^2})=(\cA_2)_{p^2q^2}
$$
where $(\cA_1)^{p^2q^2}$ is the restriction of $\cA_1$ to the factorgroup of order $p^2q^2$ whereas 
$(\cA_2)_{p^2q^2}$ is the restiction of $\cA_2$ to the subgroup of order~$p^2q^2$.
Then one can form the generalized wreath product 
$$
\cA=\cA_1\wr_{U/L}\cA_2
$$ 
over the group $\mZ_{p^2q^2r^2}$ where $U$ and $L$ are the subgroups of this group  of index and 
order~$r$ respectively. It is not diffcult to verify 
that $U/L$ is a minimal  non-quasisubprinciple $\cA$-section.\medskip

Now let $\cA$ be a quasidense S-ring over a cyclic group $G$ and $S_0$ an  $\cA$-section. Suppose 
that $\rad(\cA_{S_0})>1$  and $S_0$ is a $p$-section where $p$ is a prime divisor of $|G|$. Then 
the S-ring $\cA_{S_0}$ is the $U_0/L_0$-wreath product where $U_0$ and $L_0$ are subgroups of 
$S_0$ of index and order~$p$ respectively. Next, the sets
$$
\{U\in\cG(\cA):\ \rho_{G,S_0}(U)=U_0\}\qaq 
\{L\in\cG(\cA):\ \rho_{G,S_0}(L)=L_0\}
$$
are nonepmpty, because they contain the groups $\pi^{-1}(U_0)$ and $\pi^{-1}(L_0)$ respectively,
where $\pi=\pi_{S_0}$. Moreover, since the lattice $\cG(\cA)$ is distributive \cite[p.11]{S94}, these sets
have the greatest and least elements. Denote them by $U=U(S_0)$ and $L=L(S_0)$, respectively.
Clearly, $U\ge L$.\medskip

The definition of the section $U/L$ associated with $S_0$ is uniform.
Namely, let $S$ be an $\cA$-section  that contains $S_0$ as a subsection. Since the mapping $\rho_{G,S}$
defined in Subsection~\ref{300114a} induces a lattice epimorphism from $\cG(\cA)$ to~$\cG(\cA_S)$, we have
\qtnl{191113b}
\rho_{G,S}(U/L)=U_S/L_S
\eqtn
where $U_S/L_S$ is the section of the S-ring $\cA_S$ defined by $S_0$ in the same way as the section
$U/L$ in the S-ring~$\cA$. In particular, $U_{S_0}/L_{S_0}=U_{0^{}}/L_{0^{}}$.

\lmml{080213a}
In the above notation let $M/N\in\frS(\cA)$ be such that $M\not\le U$ and $N\not\ge L$. 
Then $S_0$ is an $\cA$-quasisubsection of~$M/N$. 
\elmm
\proof Let $S_0=V/K$. Then by the definition of $U$ and $L$ we have 
$$
U_0=(U\cap V)K/K\qaq L_0=(L\cap V)K/K.
$$ 
Since $S_0$ is a cyclic $p$-group,   this implies that 
$|V_p:U_p|=p=|L_p:K_p|$. On the other hand, from the hypothesis of the lemma it follows 
by the maximality of~$U$ that $M_p>U_p$ , and  by the minimality of~$L$ that $N_p<L_p$. Thus
$V_p\le M_p$ and $N_p\le K_p$. Therefore
$$
V_pN_{p'}\le M_p M_{p'}\le M\qaq K_pN_{p'}\ge N_pN_{p'}\ge N .
$$
Thus $T=V_pN_{p'}/K_pN_{p'}$ is a an $\cA_{M/N}$-section. Obviously, $T$ is a multiple of $V_p/K_p$. Since the latter
section is projectively equivalent to $S_0$, we are done.\bull\medskip

From Theorem~\ref{141113g} it follows that any circulant S-ring with nontrivial radical, is a proper generalized 
wreath product. Statement~(1) of the following theorem gives an explicit form of such a product in the quasidense case
(cf. Theorem~\ref{241013b}).

\thrml{291112b}
Let $\cA$ be a quasidense S-ring over a cyclic group~$G$. Suppose that $S_0\in\frS(\cA)$ is a non-quasiprincipal 
$p$-section.  Then
\nmrt
\tm{1} $\cA$ is a proper $U/L$-wreath product where $U=U(S_0)$ and $L=L(S_0)$,
\tm{2} if $T\in\frS(\cA)$ is a subsection of neither $U$ nor $G/L$, then $S_0$ is a quasisubsection of~$T$.
\enmrt
\ethrm
\proof To prove the first statement it suffices to verify that $\rad(X)\ge L$ for all $X\in\cS(\cA)_{G\setminus U}$.
Suppose on the contrary that this is not true for some~$X$. Then $S_0$ is a quasisubsection of
the section~$M/N$ where $M=\grp{X}$ and $N=\rad(X)$ (Lemma~\ref{080213a}). However, this contradicts the theorem hypothesis  because
the section $M/N$ is principal. To prove  the second statement suppose  that an $\cA$-section $T=M/N$  is  
a subsection of neither ~$U$ nor~$G/L$. Then $M\not\le U$ and $N\not\ge L$.  Thus the required statement immediately follows from Lemma~\ref{080213a}.\bull

\section{Coset S-rings}\label{051113a}

\sbsnt{Definition and basic properties.} In this section we  introduce and study circulant coset S-rings. 
In a sense these rings are antipodes of rational
circulant S-rings. Indeed, as we will see  below (Theorems~\ref{270612b} and~\ref{170611b}) an S-ring is a coset one if and only if  it can be constructed from
group rings by tensor and generalized wreath products,\footnote{In fact, the tensor product here is needless.}
whereas  an S-ring is rational if and only if  it can be constructed from S-rings of rank~$2$
in the same way (the latter follows from \cite[Theorem~1.2]{EP13}).

\dfntnl{abc}
An S-ring over an abelian group $G$ is a {\it coset} one, if any of its basic sets is a coset of a subgroup in~$G$. 
\edfntn

By definition any basic set $X$  of  a coset S-ring $\cA$ is of the form $X=xH$ for some group $H\le G$ and $x\in X$. It is easily seen that
$H=\rad(X)$, and hence
\qtnl{280513ac}
X=x\rad(X)
\eqtn
for any $x\in X$. However, $\rad(X)$ is an $\cA$-group. Thus,  any basic set of $\cA$ is a coset of a uniquely determined $\cA$-group.\medskip

The following statement expressing the ``radical monotony property" of a coset circulant S-ring, willl be used below.

\lmml{280513a}
Let $\cA$ be a circulant coset S-ring. Then given $X,Y\in\cS(\cA)$, the inclusion $Y\subset\grp{X}$ 
implies $\rad(Y)\le\rad(X)$.
\elmm
\proof Let $X,Y\in\cS(\cA)$. Then from \eqref{280513ac} it follows that $X$ and $Y$ are cosets of 
the groups $\rad(X)$ and $\rad(Y)$ respectively. Therefore the set $X^\pi$
 where $\pi=\pi_{G/\rad(X)}$, is a singleton consisting of a generator of the group $S=\grp{X}/\rad(X)$. 
This implies that $\cA_S=\mZ S$. If $Y\le\grp{X}$, then $Y^\pi$ is a basic set of $\cA_S$, and hence $Y^\pi\subset S$ is also a singleton. It follows that 
$\rad(Y)^\pi\le \rad(Y)^\pi=1$. Thus $\rad(Y)\le\rad(X)$ as required.\bull\medskip

The circulant coset S-rings can be characterized in terms of their sections as follows.

\thrml{270612b}
For a circulant S-ring $\cA$ the following statements are equivalent:
\nmrt
\tm{1} $\cA$ is a coset S-ring,
\tm{2} $\cA_S=\mZ S$ for any principal $\cA$-section $S$,
\tm{3} $\cA_S=\mZ S$ for any $\cA$-section $S$ with trivial radical,
\tm{4} $\cA_S=\mZ S$ for any quasisubprincipal $\cA$-section $S$.
\enmrt
\ethrm
\proof Statements (1) and (2) are equivalent: implication $(1)\Rightarrow(2)$ follows from~\eqref{280513ac} whereas implication $(2)\Rightarrow(1)$
follows from the definition of principal section. Next, implication $(4)\Rightarrow(2)$ is obvious and implication $(2)\Rightarrow(4)$ is true because
the equality in statement~(4) is preserved under  projective equivalence and taking subsections. Finally, 
any principal $\cA$-section obviously 
has trivial radical, and any $\cA$-section with trivial radical is subprincipal (Lemma~\ref{100214a}). Thus  
implications $(3)\Rightarrow(2)$ and $(4)\Rightarrow(3)$ hold.\bull\medskip

Any primitive section $S$ of a dense circulant S-ring $\cA$ has prime order. Therefore $\rad(\cA_S)=1$. Thus
if $\cA$  is a coset S-ring, then $\cA_S=\mZ S$ for any primitive $\cA$-section $S$ (Theorem~\ref{270612b}). The converse statement is not true,  a counterexample is given 
by $\cA=\cyc(\{\pm 1\},\mZ_8)$.

\thrml{170611b}
The class of circulant coset S-rings is closed under  restriction to a section and under tensor and generalized wreath products, 
and consists of quasidense S-rings.
\ethrm
\proof Since any quotient epimorphism takes a coset to a coset, and the product of cosets is a coset, the closedeness  statement follows
from the definitions of tensor and generalized wreath products. The quasidensity 
statement is true because any non-quasidense S-ring has a rank~$2$ section of composite order 
(Theorem~\ref{060214a}) whereas by
above no coset S-ring can have such a section.\bull\medskip

The intersection of circulant coset S-rings is not necessarily a coset one: a counterexample is given 
by the S-ring $\cA=\cA_1\cap\cA_2$  over the group $\mZ_{pq}$ where 
$$
\cA_1=\mZ\mZ_p\wr\mZ\mZ_q\qaq \cA_2=\mZ\mZ_q\wr\mZ\mZ_p
$$
with $p$ and $q$ distinct primes. One can see that its rank equals~$2$, and
hence it is not coset. Moreover, $\cA$ is even not quasidense (Theorem~\ref{060214a}). The following statement shows 
that this is a unique obstacle.

\thrml{311013a}
The intersection of circulant coset S-rings is a coset one whenever it contains a quasidense S-ring. 
\ethrm
\proof Let $\cA=\cA_1\cap\cA_2$ where $\cA_1$ and $\cA_2$ are circulant coset S-rings. Suppose that $\cA$
contains a quasidense S-ring. Then by Theorem~\ref{180214a} we can assume that $\cA$ is quasidense.
By implication
$(3)\Rightarrow(1)$ of Theorem~\ref{270612b} it suffices to verify that $\cA_S=\mZ S$ for an $\cA$-section $S$ 
with trivial radical. However, since $\cA_i\ge\cA$ for $i=1,2$, any such section is an $\cA_i$-section and 
$(\cA_i)_S\ge\cA_S$.
The quasidensity of the S-ring $\cA_S$ implies by Theorem~\ref{060214b} that $\rad((\cA_i)_S)=1$. 
Taking into account that $\cA_i$ is a coset S-ring, we conclude that $(\cA_i)_S=\mZ S$ (implication $(1)\Rightarrow(3)$ 
of Theorem~\ref{270612b}). Thus by~\eqref{311013b} we have
$$
\cA_S=(\cA_1)_S\cap(\cA_2)_S=\mZ S
$$
as required.\bull\medskip

In what follows any representation of an S-ring $\cA$ as a proper generalized wreath product, will be called 
a {\it gwr-decomposition} 
of~$\cA$. Given an $\cA$-section $S$ we say that the $T'$-decomposition of $\cA_S$ {\it is lifted} to
a $T$-decomposition of~$\cA$ if $\cA$ is the $T$-wreath product and $T_S=T'$ (cf.~Theorem~\ref{031212c}).

\thrml{290513a}
Let  $S$ be a section of a circulant coset S-ring $\cA$. Then any  gwr-decomposition of the S-ring~$\cA_S$ can be lifted to
a gwr-decomposition of~$\cA$.
\ethrm
\proof Let $\cA_S$ be a proper $V/K$-wreath product. To lift it to a gwr-decomposition of~$\cA$ it suffices to consider two cases
depending on whether $S$ is  a subgroup or quotient of~$G$. By duality  the latter case follows from the former one by
statement~(2) of Lemma~\ref{310114a}
and statement~(3) of Theorem~\ref{140509a}. In the former case $V$ and $K$ are also subgroups of~$G$. Set  $L=K$ and
$$
U=\lg\{X\in\cS(\cA):\ \lg X\rg\cap S\le V\}\rg.
$$
Clearly, $U\in\cG(\cA)$ and $U\ge V\ge L$. On the other hand, if $H_1$ and $H_2$ are subgroups of $G$ such that 
$H_1\cap S\le V$ and $H_2\cap S\le V$, then by the  distributivity of the lattice~$\cG(\cA)$ we have
$$
H_1H_2\cap S=(H_1\cap S)(H_2\cap S)\le V.
$$
This shows that $U\cap S\le V$, and hence $U\cap S=V$. In particular, $U\ne G$. To complete the proof
let $X$ be a basic set of $\cA$ outside $U$. Then by the definition of~$U$ the group $H=\lg X\rg\cap S$ 
is not a subgroup of~$V$. Since $\cA_S$ is the $V/K$-wreath product, this implies that $K\le\rad(Y)$ 
where $Y$ is a highest basic set of~$\cA_H$. However, then by Lemma~\ref{280513a} we have 
$$
L=K\le\rad(Y)\le\rad(X).
$$ 
Thus $\cA$ is a proper $U/L$-wreath product. Since also $\rho_{G,S}(U/L)=V/K$, we are done.\bull

\sbsnt{Elementary coset S-rings.}
From Theorems~\ref{270612b} and~\ref{170611b} it follows that any circulant coset S-ring can be constructed from
group rings by generalized wreath products. In the rest of this section we are interested in the coset S-rings that
are obtained in  one iteration of the above process. More precisely,  by Theorem~\ref{230114a} given a section $T=U/L$ 
of a cyclic group~$G$ one can form the S-ring
\qtnl{311013e}
\mZ(G,T)=\mZ U\wr_T\mZ(G/L)
\eqtn
because any group ring is dense and the restrictions of both $\mZ U$ and $\mZ(G/L)$ to $T$ equal $\mZ T$.
It is easily seen that  $\mZ(G,T)$ is a coset S-ring over~$G$. 

\dfntn
Any S-ring of the form  \eqref{311013e} is called elementary coset.
\edfntn

Clearly, the group ring $\mZ G$ is elementary coset (in this case generalized wreath product~\eqref{311013e} is not 
proper). It is also easily seen that every basic set of elementary coset S-ring~\eqref{311013e} inside $U$ is a
singleton whereas  the basic sets outside $U$ are $L$-cosets. 
Any elementary coset S-ring is schurian (Theorem~\ref{060813a}), and the automorphism group of it 
is the canonical generalized wreath product of $U_{right}$ by $(G/L)_{right}$ in the sense of~\cite{EP12}. 
For association schemes a similar situation was studied in~\cite{M09}.\medskip

Let $\cA$ be elementary coset S-ring~\eqref{311013e}. Given a function $t\in L^{G/U}$ and an element $g\in G$ 
we define a permutation $\sigma_{t,g}\in\sym(G)$ by
\qtnl{170613d}
\sigma_{t,g}:x\mapsto xt(X)g,\quad x\in G,
\eqtn 
where $X$ is the $U$-coset  containing~$x$. The set of all these permutations forms a group with identity
$\sigma_{1,1}$ where the first $1$ in subscript denotes the function taking every~$x$ to~$1$, and the multiplication
satisfying
\qtnl{250314a}
\sigma_{t_1,g_1}\sigma_{t_2,g_2}=\sigma_{t_1t_2,g_1g_2}
\eqtn
for all $t_1,t_2\in L^{G/U}$ and $g_1,g_2\in G$. Clearly, $\{\sigma_{1,g}:\ g\in G\}=G_{right}$ and 
$\sigma_{1,g}=\sigma_{t,1}$ for all $g\in L$ where $t$ is the function taking every~$x$ to~$g$.\medskip

It follows from~\cite[Theorem~7.2]{EP10} that 
$\sigma_{t,g}\in\aut(\cA)$ whenever $t(U)=g=1$.  Since $G_{right}\le\aut(\cA)$,  the permutation 
$\sigma_{t,g}$ is an automorphism of~$\cA$
for all $t$ and $g$. In fact, statement~(2) of the theorem below shows that $\cA$ has no other automorphisms.

\thrml{300413y}
Let $\cA=\mZ(G,T)$ be an elementary coset S-ring~\eqref{311013e} and $S$ an $\cA$-section. Then
\nmrt
\tm{1} $\cA_S=\mZ(S,T_S)$ where $T_S$ is the section defined in~\eqref{170613c},
\tm{2} $\aut(\cA)=\{\sigma_{t,g}:\ t\in L^{G/U},\ g\in G\}$,
\tm{3}  $\aut(\cA)^S=\aut(\cA_S)$.
\enmrt
\ethrm
\proof By Theorem~\ref{031212c} the S-ring $\cA_S$ is the $T_S$-wreath product
of the S-rings $\cA_{U_S}$ and $\cA_{S/L_S}=\cA_{(G/L)_S}$ which are Cayley isomorphic 
by Lemma~\ref{210613a} and Theorem~\ref{261010a}
to the S-rings $\cA_{S_U}$ and $\cA_{S_{G/L}}$
respectively. Since $\cA_U=\mZ U$ and $S_U$ is an $\cA_U$-section, as well as $\cA_{G/L}=\mZ(G/L)$ and $S_{G/L}$ is an
$\cA_{G/L}$-section, we have  $\cA_{S_U}=\mZ S_U$  and $\cA_{S_{G/L}}=\mZ (S_{G/L})$ which proves statement~(1).\medskip

To prove statement~(2) denote by $\Gamma$ the group in the right-hand side of the equality, Then from the discussion before the theorem
it follows that $\Gamma\le\aut(\cA)$. 
Therefore to check the reverse inclusion it suffices to prove that $|\aut(\cA)|\le|\Gamma|$. However, since $\cA_{G/L}=\mZ (G/L)$,
we have
$$
\aut(\cA)^{G/L}=(G/L)_{right}=\Gamma^{G/L}.
$$
To complete the proof we show that the kernel of the epimorphism 
$$
\pi:\aut(\cA)\to\aut(\cA)^{G/L}
$$ 
is contained in~$\Gamma$,
more precisely that any $\sigma\in\ker(\pi)$ is of the form $\sigma_{t,1}$ for some $t\in L^{G/U}$.  Note
that such a permutaton~$\sigma$ leaves each $X\in G/U$ fixed as a set. Therefore all we need to prove is that
$\sigma^X$ acts on $X$ by multiplying by an element $l=l(X)$ belonging to~$L$. Let $g_X\in G_{U\to X}$. Then obviously the permutation  $\sigma'=g_X^{}\sigma g_X^{-1}$ belongs to $\ker(\pi)$ and leaves the set
$U$ fixed. Taking into account that $\cA_U=\mZ U$, we see that $(\sigma')^U\in\aut(\cA)^U=U_{right}$. 
Therefore since $\sigma'$ leaves also any $L$-coset fixed,
there exists $l\in L$ for which
$$
u^{g^{}_X\sigma g_X^{-1}}=u^{\sigma'}=lu,\qquad u\in U.
$$
When $u$ runs over $U$, the element $x=u^{g_X}$ runs over $X$. Therefore the above equality implies that 
$x^\sigma=lx$ for all $x\in X$, as required.\medskip

To prove statement~(3) we observe that obviously $\aut(\cA)^S\le\aut(\cA_S)$. To verify the converse
inclusion, let $f'\in \aut(\cA_S)$ and $S=H/K$. Then by statements~(1) and~(2) we have $f'=\sigma_{t',g'}$
where $t'\in L_S^{S/U_S}$ and $g'\in S$. However, from~\eqref{250314a} it follows that 
$\sigma_{t',g'}=\sigma_{t',1}\sigma_{1,g'}$. Moreover, the permutation $\sigma_{1,g'}\in S_{right}$
can be lifted to a permutation in $G_{right}$, because obviously $\Gamma^S=S_{right}$ where $\Gamma$
is the setwise stabilizer of~$H$ in~$G_{right}$. Thus  we can assume that $f'=\sigma_{t',1}$.\medskip

To verify that $f'$ can be lifted to a permutation in $\aut(\cA)$ it suffices to check  that there exists  
$t\in L^{G/U}$ such that
\qtnl{240314a}
(\sigma_{t,1})^S=\sigma_{t',1}.
\eqtn
To define the function~$t$ given
$X\in G/U$ set $t(X)=1$ if the set $X\cap H$ is empty. Suppose that it is not empty. Then
$X\cap H$ is a coset of $H\cap U$. Recall that
$$
U_S=(U\cap H)K/K\qaq L_S=(L\cap H)K/K.
$$
Therefore, the first equality  implies that $X':=\pi(X\cap H)$ belongs to~$S/U_S$ where $\pi=\pi_S$.  
Moreover, the set $L\cap \pi^{-1}(t'(X')$ is not empty by the second equality. Now set $t(X)$ to be any element of that set.
It follows that in any case $t\in L^{G/U}$, and  the permutation $\sigma_{t,1}$ leaves $H$ fixed. Thus,
given $x'\in S$ we have
$$
(x')^{\sigma_{t,1}^S}=x'\pi(t(X))=x't'(X')=(x')^{\sigma_{t',1}}
$$
where $X'$ is the $U_S$-coset containing~$x'$, which proves~\eqref{240314a}.\bull\medskip

In general, an extension of a coset S-ring is not necessary coset: a non-coset S-ring $\cyc(\{\pm 1\},\mZ_8)$
contains a coset S-ring $\cyc(K,\mZ_8)$ where $K=\aut(\mZ_8)$.  However, we can prove the 
following auxiliary statement to be used in Section~\ref{200214b}.

\lmml{280314a}
Let $\cA$ be an elementary coset S-ring~\eqref{311013e}. Suppose that either $|L|$ or $|G/U|$ is prime. Then 
any extension of~$\cA$ is an elementary coset S-ring.
\elmm
\proof By duality (Theorem~\ref{140509a}) without loss of generality we can assume that $|L|$ is prime.
Let $\cA'\ge\cA$. Denote by $U'$ the group consisting of all $x\in G$ for which $\{x\}$ is a basic set of~$\cA'$.
Clearly, $U'$ is an $\cA'$-group containing~$U$. We claim that any basic set not in $U'$, is a basic set of $\cA'$. 
Then since $\cA'\ge\cA$, we have $\cA'=\mZ(G,T')$ where $T'=U'/L$, as required.\medskip

Suppose that the claim  is not true. Then there exists a basic set $X\not\subset U'$ of~$\cA$, that is not
a basic set of $\cA'$. It follows that $X\subset G\setminus U$. Since $\cA=\mZ(T,G)$, we have
\qtnl{280314s}
X=xL
\eqtn
for all $x\in X$. Moreover, $X$   is the union of at least two basic sets of $\cA'$. If one of them, say $X'$, has nontrivial 
radical, then the latter coincides with $L$, because $\cA'_{G/L}\ge\cA^{}_{G/L}=\mZ(G/L)$. But then 
from~\eqref{280314s} it follows that
$$
|L|=|X|\ge |X'|\ge|L|.
$$
Therefore $X=X'$. Contradiction.\medskip

Now we have $\rad(X')=1$ for all $X'\in\cS(\cA')_X$. Since the group $H:=\grp{X}$ is cyclic, there exists such an
$X'$ such that $\grp{X'}=H$. Therefore $\rad(\cA'_{H})=1$. On the other hand, since $\cA^{}_H\le\cA'_H$
and $\cA$ is quasidence, we conclude by Theorem~\ref{180214a}, then $\cA'_H$ is also quasidence.
Thus by Theorem~\ref{060214b} we have 
\qtnl{260314e}
\cA'_H=\cyc(K,H)
\eqtn
for some group $K\le\aut(H)$. Clearly, $K$ is a subgroup of the stabilizer of~$1$ in $\aut(\cA_H)$. Since 
this stabilizer is a $p$-group (statements~(1) and~(2) of Theorem~\ref{300413y}) where $p=|L|$, it follows 
that $K$ is also a $p$-group. This implies that the cardinality of every $K$-orbit is a $p$th power. 
Since $X$ is a $K$-invariant set of cardinality $p$, equality~\eqref{260314e} implies that
either every element in $\cS(\cA')_X$ is singleton, or $\cS(\cA')_X=\{X\}$. In both cases we come to a
contradiction with the choice of~$X$.\bull

\section{Schurity and separability of coset S-rings}\label{200214b}

The theory of coset circulant S-rings developed in Section~\ref{051113a} enables us to prove the following theorem which
is the main result of the section.

\thrml{150113a}
Any circulant coset S-ring is schurian and separable.
\ethrm

We will deduce Theorem~\ref{150113a} in the end of the section from two following auxiliary statements 
on the automorphism group of a coset S-ring. In what follows, if an S-ring~$\cA$ is the $S$-wreath product
where $S\in\frS(\cA)$, then we say that $S$ is a {\it gwr-section} of~$\cA$.

\thrml{300413w}
Let $\cA$ be a coset S-ring over a cyclic group $G$. Then the group $\aut(\cA)$ is generated by the automorphism groups 
of elementary coset S-rings
$\mZ(G,S)$ where $S$ runs over all   gwr-sections of $\cA$.
\ethrm
\proof Induction on the order of $G$. The base case of the induction is obvious. Let $G\ne 1$ and $f\in\aut(\cA)$. By 
Theorem~\ref{170611b} the S-ring $\cA$ is quasidense. So
there exists an $\cA$-group $L$ of prime order. Set $S=G/L$. Since $f^S\in\aut(\cA_S)$, by the inductive hypothesis 
$f^S$ is the product of automorphisms 
of elementary coset S-rings $\mZ(S,T')$ where $T'$ runs over all   gwr-sections of $\cA_S$. On the other hand, by Theorem~\ref{290513a} any
$T'$-decomposition of $\cA_S$ can be lifted to a $T$-decomposition of~$\cA$.
Set $\cB=\mZ(G,T)$. Then $\cB_S=\mZ(S,T')$ by statement~(1) of Theorem~\ref{300413y}.
Moreover, by statement~(3) of that theorem we also have $\aut(\cB)^S=\aut(\cB_S)$. Thus  to write $f$ as the product
of automorphisms of elementary coset S-rings, without loss of generality we can assume that
\qtnl{110613a}
f^S=\id_S.
\eqtn

Denote by $\cA'$ the S-ring over $G$ associated with the group $\Gamma=\grp{G_{right},f}$. 
In particular, $f\in\aut(\cA')$. Due to~\eqref{110613a} we have $\Gamma^S=S_{right}$. Therefore
\qtnl{270314a}
\cA'_{G/L}=\mZ(G/L).
\eqtn
Moreover,  $\cA'\ge\cA$ because $f\in\aut(\cA)$. It follows that  $\cA'_L\ge\cA^{}_L$. On the other hand, 
$\cA_L$ is a coset S-ring (Theorem~\ref{170611b}) with trivial radical because $|L|$ is prime. 
Thus $\cA_L=\mZ L$, and hence  $\cA'_L=\mZ L$.  Together with~\eqref{270314a} this implies that $\cA'$ 
is an extension of elementary
coset S-ring $\mZ(G,L/L)$. Therefore by Lemma~\ref{280314a} we conclude that
\qtnl{110613b}
\cA'=\mZ(G,T')
\eqtn
where $T'=U'/L$ with $U'$ being the largest  $\cA'$-group for which $\cA'_{U'}=\mZ U'$. \medskip

Denote by $U$ the largest $\cA$-subgroup of $G$ inside $U'$. Clearly, $L\le U$.  Take a set $X\in\cS(\cA)$ 
outside~$U$. Then by the definition of $U$ there exists a basic $\cA'$-set $X'\subset X$ that is outside $U'$. 
By~\eqref{110613b} this implies that $X'$ is an $L$-coset. However, $X$ is a coset of $\rad(X)$
because $\cA$ is a coset S-ring. Thus $L\le\rad(X)$. This proves that
$T:=U/L$ is a gwr-section of~$\cA$. Besides, since $U\le U'$
we have
$$
\mZ(G,T)\le\mZ(G,T').
$$
Since $f\in\aut(\cA')$, this inclusion together with \eqref{110613b} implies that
$f$ is an automorphism of the elementary coset S-ring $\mZ(G,T)$, as required.\bull

\thrml{300413v}
Let $\cA$ be a circulant coset S-ring. Then $\aut(\cA)^S=\aut(\cA_S)$ for any $\cA$-section~$S$.
\ethrm
\proof Obviously, $\aut(\cA)^S\le\aut(\cA_S)$. To prove the reverse inclusion we observe that by Theorem~\ref{300413w}
each automorphism of $\cA_S$  can be written as the product of automorphisms of elementary coset S-rings $\mZ(S,T')$ 
where $T'$ runs over all   gwr-sections of $\cA_S$. On the other hand, by Theorem~\ref{290513a} any $T'$-decomposition 
of $\cA_S$ can be lifted to a $T$-decomposition of~$\cA$. Thus the required statement follows from 
statement~(3) of Theorem~\ref{300413y}.\bull\medskip

{\bf Proof of Theorem~\ref{150113a}.} Let $\cA$ be a coset S-ring over a cyclic group~$G$. We will prove both statements 
by induction on the order of this group. Let $|G|>1$.
If $\rad(\cA)=1$, then $\cA=\mZ G$ by implication $(1)\Rightarrow(3)$ of 
Theorem~\ref{270612b} and the statements are obvious. Let now $\rad(\cA)>1$. Then 
by statement~(1) of Theorem~\ref{141113g} there is a  gwr-decomposition
\qtnl{150513a}
\cA=\cA_U\wr_{U/L}\cA_{G/L}.
\eqtn 
By Theorem~\ref{170611b} the  S-rings $\cA_U$ and $\cA_{G/L}$ are coset ones. Thus they are
schurian and separable by the inductive hypothesis.\medskip

To prove that $\cA$ is schurian set $\Delta_1=\aut(\cA_U)$ and $\Delta_0=\aut(\cA_{G/L})$. 
Then by Theorem~\ref{060813a} it suffices to verify that 
$(\Delta_1)^{U/L}=(\Delta_0)^{U/L}$. However, from Theorem~\ref{300413v} it follows that
$$
(\Delta_1)^{U/L}=\aut(\cA_U)^{U/L}=\aut(\cA_{U/L})=\aut(\cA_{G/L})^{U/L} =(\Delta_0)^{U/L}
$$
as required.\medskip

To prove that $\cA$ is separable let $\varphi:\cA\to\cA'$ be a similarity. Then by statement~(1) of Theorem~\ref{271213x}
the S-ring $\cA'$ is the $U'/L'$-wreath product where $U'=U^\varphi$ and $L'=L^\varphi$. By the inductive hypothesis
the similarities 
$$
\varphi_U:\cA^{}_{U^{}}\to\cA'_{U'}\qaq\varphi_{G/L}:\cA^{}_{G^{}/L^{}}\to\cA'_{G'/L'}
$$ 
are induced by some bijections, say $f_1$ and $f_0$. Given $X\in G/U$ let us fix two bijections 
$g\in G^{}_{U^{}\to X^{}}$ and 
$g'\in G'_{X'\to U'}$.
By Lemma~\ref{220513a} the bijection $g^{U/L}(f_0)^{X/L}(g')^{X'/L'}$ induces the similarity $\varphi^{U/L}$. It follows
that
\qtnl{220513z}
g^{U/L}(f_0)^{X/L}(g')^{X'/L'}((f_1)^{U/L})^{-1}\in\aut(\cA_{U/L}).
\eqtn
By Theorem~\ref{300413v} there exists an automorphism $h_X\in\aut(\cA_U)$ such that its restriction to $U/L$ coincides with the automorphism in
the left-hand side of~\eqref{220513z}. Now, the family of bijections
$$
f_X:=g^{-1}h_Xf_1(g')^{-1},\qquad X\in G/U,
$$
defines a bijection $f:G\to G'$ such that $(G/U)^f=G'/U'$ and $(G/L)^f=G'/L'$. Then obviously 
\qtnl{310314a}
gf^Xg'=gf_Xg'=h_Xf_1.
\eqtn
By~\eqref{150114a} this implies that the bijection $gf^Xg'$
induces the similarity $\varphi_U$ for all $X$. This proves the second part of condition~\eqref{250613b}. 
Let us check the first part of that condition. For any $X\in G/U$ from the definition of $h_X$ and~\eqref{310314a}
it follows that
$$
f^{X/L}=(g^{-1}h_Xf_1(g')^{-1})^{X/L}=
$$
$$
(g^{U/L})^{-1}g^{U/L}(f_0)^{X/L}(g')^{X'/L'}((f_1)^{U/L})^{-1}(f_1)^{U/L}((g')^{X'/L'})^{-1}=(f_0)^{X/L}.
$$
Therefore $f^{G/L}=f_0$, and hence the bijection $f^{G/L}$ induces the similarity $\varphi_{G/L}$. Thus by 
Theorem~\ref{140513a} the bijection $f$ induces the similarity $\varphi$ as required.\bull\medskip

The following statement can in a sense be regarded as a combinatorial analog of  Theorem~\ref{150113a}.

\crllrl{120513s}
Any circulant coset S-ring is an intersection of elementary coset S-rings over the same group.
\ecrllr
\proof Let $\cA$ be a circulant coset S-ring. Denote by $\cA'$ the intersection of all elementary coset S-rings
$\mZ(G,S)$ where $S$ runs over all   gwr-sections of $\cA$.  Since $\cA$ is quasidense (Theorem~\ref{170611b}),
the S-ring~$\cA'$ is coset (Theorem~\ref{311013a}). So
from Theorem~\ref{150113a}  it follows that
it is schurian. Therefore~$\cA'$ equals the S-ring associated with the group $\Gamma$ generated 
by the automorphism groups 
of the above S-rings $\mZ(G,S)$. On the other hand,  by Theorem~\ref{150113a}  the S-ring $\cA$ is also schurian. By 
Theorem~\ref{300413w} this implies that $\cA$ equals the S-ring associated with the group $\Gamma$. 
Thus $\cA=\cA'$.\bull

\section{Coset closure}\label{111113b}

\sbsnt{Relative coset closure.}
We start with developing a technique to find the schurian closure of a quasidense circulant S-ring~$\cA$
(Theorem~\ref{251212a}). The key point of this technique is its coset closure $\cA_0$
defined in Introduction (Definition~\ref{201113a}) as the intersection of all coset S-rings containing~$\cA$. By 
Theorem~\ref{311013a} the S-ring~$\cA_0$ is a coset one. Moreover,
\qtnl{201113b}
(\cA_0)_S=\mZ S,\qquad S\in\frS_0(\cA)
\eqtn
where $\frS_0(\cA)$ is the class of all quasisubprincipal $\cA$-sections.  Indeed, by Theorem~\ref{261010a} this
is reduced to the case of a principal $S$, in which the required statement follows
from Theorem~\ref{060214b} and implication $(1)\Rightarrow(3)$ of Theorem~\ref{270612b}.
Thus the S-ring~$\cA_0$ equals the intersection of all coset S-rings~$\cA'$ such that
\qtnl{141113a}
\cA'\ge\cA \qaq (\cA')_S=\mZ S
\eqtn 
for all $S\in\frS_0(\cA)$. For induction reasons it is convenient to generalize the coset closure concept 
by permitting the section $S$  to run over a larger class~$\frS$.

\dfntn
A class $\frS\subset\cF(\cA)$ is {\it admissible} with respect to~$\cA$ (or $\cA$-admissible), 
if $\cF_{prin}(\cA)\subset\frS\subset\cF_{cyc}(\cA)$
and $\frS$ is closed under taking $\cA$-quasisub\-sections.
\edfntn

Given an $\cA$-admissible class $\frS$ and an $\cA$-section $S$ we set $\frS_S=\rho_{G,S}(\frS)$
where the mapping~$\rho_{G,S}$ is defined in Subsection~\ref{300114a}.  Since $\rho_{G,S}$
is identical on the set $\frS(S)\subset\frS(G)$, we have
\qtnl{281113a}
\frS_S=\{\rho_{G,S}(S'):\ S'\in\frS,\ S'\preceq S\}.
\eqtn
The following statement gives the properties of admissible classes.

\lmml{141113v}
Let $\frS$ be an admissible class with respect to a quasidense circulant S-ring~$\cA$. Then
\nmrt
\tm{1} $\frS_0(\cA)\subset\frS$,
\tm{2} the class $\frS_S$ is $\cA_S$-admissible for any $S\in\cF(\cA)$,
\tm{3} if $S\in\cF(\cA)\setminus\frS$, then $\cA_S$ admits a gwr-decomposition.
\enmrt
\elmm
\proof Statements (1) is obvious. Statement (2)  follows from \eqref{281113a}. To prove Statement~(3) 
let $S\in\cF(\cA)\setminus\frS$. Then $\rad(\cA_S)\ne 1$, for otherwise $S$ is a subsection of a principal 
$\cA$-section (Lemma~\ref{100214a}) whereas $\frS$
contains all such sections. Thus the required statement follows from Theorem~\ref{141113g}.\bull\medskip

Let $\cA$ be a quasidence circulant S-ring and $\frS$ an $\cA$-admissible class. Then there is at least
one coset S-ring $\cA'$, namely the group ring, for which relations~\eqref{141113a} hold
for all $S\in\frS$. This justifies the following definition.

\dfntn
The {\it coset closure} $\cA_{0,\frS}$  of $\cA$ with respect to  $\frS$
is the intersection of all coset S-rings $\cA'$ such that relations~\eqref{141113a}
hold for all $S\in\frS$. 
\edfntn

Clearly, $\cA_{0,\frS}\ge\cA$. Moreover, from~\eqref{311013b} it follows that $(\cA_{0,\frS})_S=\mZ S$ for
all $S\in\frS$. Besides, the discussion in the first paragraph of the section shows that
\qtnl{191113s}
\cA_0=\cA_{0,\frS}\quad\text{when}\quad \frS=\frS_0(\cA).
\eqtn

The sense of the following definition will be clarified in Corollary~\ref{221113a}. 
With any $\cA$-admissible class $\frS$  we associate a larger class of sections
defined as follows
$$
\frSH=\{S\in\frS(\cA):\ S_p\in\frS\ \text{for all primes}\ p\ \text{dividing}\ |S|\}. 
$$
The class $\frSH$ is closed under taking quasisubsections, but generally is not admissible because may
contain non-cyclotomic sections. As the example on page~\pageref{230414a} shows,
in general it can be larger than $\frS$.


\lmml{031212a}
Let $\cA$ be a quasidence S-ring over a cyclic group~$G$.  Then for any $\cA$-admissible class $\frS$ we have
\nmrt
\tm{1} the S-ring $\cA_{0,\frS}$ is a coset one,
\tm{2} $(\cA_{0,\frS})_S=\mZ S$ for all $S\in\frS$,
\tm{3} any coset belonging to a section in $\frSH$ is an $\cA_{0,\frS}$-set.
\enmrt
\elmm
\proof Statement~(1) immediately follows from Theorem~\ref{311013a}. Statement~(2) follows from
the remark after the definition of the coset closure $\cA_{0,\frS}$.
To prove statement~(3) let $\cA'$ be a coset S-ring such that
relations~\eqref{141113a} hold for all $S\in\frS$. Then given a section $S\in\frSH$ we have 
$S_p\in\frS$, and hence $(\cA')_{S_p}=\mZ S_p$
for all primes $p$ dividing $|S|$. This implies that $(\cA')_S=\mZ S$. It follows that if $S=U/L$, then any
$L$-coset in $U$ is an $\cA'$-set. Since $\cA$ is the intersection of all such $\cA'$,  it is an $\cA$-set
as required.\bull

\sbsnt{Lifting.} Let $\frS$ be an admissible class with respect to a quasidence S-ring $\cA$  over a cyclic group $G$. 
Suppose that 
$\cA$ is the $S$-wreath product where $S=U/L$ is an $\cA$-section. We say that this product
is {\it $\frS$-consistent} if any section in $\frS$ is either an $\cA_U$- or $\cA_{G/L}$-section. 
Below to simplify notation we write $(\cA_S)_{0,\frS}$ instead of $(\cA_S)_{0,\frS_S}$.

\thrml{291112u}
Let $\cA$ be a circulant quasidense S-ring, $\frS$ an $\cA$-admissible class and $S$ an $\cA$-section.
Then the following conditions are equivalent:
 \nmrt
\tm{1} $(\cA_S)_{0,\frS}\ne \mZ S$,
\tm{2} $S\not\in\frSH$,
\tm{3} there exists an $\frS_S$-consistent  gwr-decomposition of $\cA_S$,
\tm{4} there exists an $\frS_S$-consistent  gwr-decomposition of $\cA_S$ that can be lifted to
an $\frS$-consistent gwr-decomposition of~$\cA$.
\enmrt
\ethrm
\proof Let us prove implication $(1)\Rightarrow(2)$. Suppose on the contrary that $S\in\frSH$.
Then obviously $S\in\wh{\frS_S}$. By statement~(3) of Lemma~\ref{031212a} applied to $\cA=\cA_S$,
$\frS=\frS_S$ and the section $S/1$, this implies that $(\cA_S)_{0,\frS}=\mZ S$. Contradiction.\medskip

To prove implication $(2)\Rightarrow(4)$ let $S\not\in\frSH$. Then there is  a prime divisor~$p$ of~$|S|$ such that 
$S_p\not\in\frS$.  Suppose first
that $S_p$ is not an $\cA_S$-group. Then the S-rings $\cA_S$ and $\cA$ are not dense. Then by Theorem~\ref{301012a}  there
exists an elementary  nondense $\cA_S$-section~$S_0$. By the first part of Theorem~\ref{290113a} applied to
 the S-rings $\cA$ and $\cA_S$, and equality~\eqref{191113a} there are gwr-decompositions
\qtnl{050413a}
\cA=\cA_U\wr_{U/L}\cA_{G/L}\qaq \cA_S=\cA_{U_S}\wr_{U_S/L_S}\cA_{S/L_S}
\eqtn
where $U=U(S_0)$ and $L=L(S_0)$ and $U_S/L_S=\rho_{G,S}(U/L)$. It follows that the first one is a lifting of the second.
Finally, the second part of Theorem~\ref{290113a} together with the fact that
any section in the class $\frS$ is dense, shows that these gwr-decompositions are 
 $\frS$- and $\frS_S$-consistent respectively.
Thus statement~(4) holds in this case.\medskip

Let now $S_p$ be  an $\cA_S$-group. Then by Lemma~\ref{141113v} the hypothesis of Theorem~\ref{291112b} holds for the S-rings $\cA$ and $\cA_S$
with $S_0=S_p$ in both cases. So by statement~(1) of this theorem there are gwr-decompositions~\eqref{050413a} and 
due to~\eqref{191113b} the first one is a
lifting of the second. To prove that the first decomposition is  $\frS$-consistent, suppose on the contrary that there exists
$T\in\frS$ which is neither $\cA_U$- nor $\cA_{G/L}$-section. Then by  statement~(2) of that theorem 
$S_p$ is a quasisubsection of~$T$. Therefore
$S_p\in\frS$,  which contradicts the assumption on~$S_p$. The $\frS$-consistency of the second decomposition is proved
similarly. Thus statement~(4) holds in this case too.\medskip

Implication $(4)\Rightarrow(3)$ is obvious. To prove implication $(3)\Rightarrow(1)$ without loss of generality we
can assume that $S=G$. Suppose that the S-ring $\cA$ admits an $\frS$-consistent $U/L$-decomposition. 
By Theorem~\ref{230114a} we can form the S-ring
$$
\cB=\cA'_U\wr^{}_{U/L}\cA'_{G/L}
$$
where $\cA'=\cA_{0,\frS}$. By Theorem~\ref{170611b} this S-ring is a coset one. Moreover, the consistency property 
implies that any section $T\in\frS$ is either an
$\cA_U$- or $\cA_{G/L}$-section, and hence either an $\cA'_U$- or $\cA'_{G/L}$-section. Therefore 
$\cB_T=\mZ T$ for all such~$T$. Thus, by the definition of the 
coset closure  we have $\cB\ge \cA'$. So from  the minimality property of the generalized wreath product 
it follows that the S-ring $\cA'$ admits the $U/L$-decomposition. Therefore
$\cA'\ne \mZ G$ as required.\bull\medskip

The following auxiliary statement will be used in proving the  theorems in the next subsection.

\lmml{241212a}
Let $\cA$ be a circulant quasidense S-ring, $\frS$ an $\cA$-admissible class and $S$ an $\cA$-section. Suppose that
$\cA_S$  admits an $\frS_S$-consistent $T$-decom\-position that can be lifted to an $\frS$-consistent gwr-decomposition
of~$\cA$. Then $(\cA_{0,\frS})_S$ admits the $T$-decomposition. 
\elmm
\proof By Theorem~\ref{031212c} without loss of generality we can assume that $S=G$. Then we have to verify
$\cA_{0,\frS}$ admits the $T$-decomposition whenever $\cA$ admits an $\frS$-consistent $T$-decom\-position.
However, the latter implies that any $T'\in\frS$ is either an
$\cA_U$- or  $\cA_{G/L}$-section where $U/L=T$. It follows that $T'$ is either $\cA'_U$- or $\cA'_{G/L}$-section
where $\cA'$ is the elementary coset S-ring $\mZ(G,T)$. Therefore
$$
(\cA')_{T'}=\mZ {T'},\quad T'\in\frS.
$$
So by the definition of the coset closure we have $\cA'\ge\cA_{0,\frS}$. Since obviously
$T$ is  an $\cA_{0,\frS}$-section, the S-ring $\cA_{0,\frS}$ is the $T$-wreath product, as required.\bull

\sbsnt{Main properties.}
Here we  study the coset closure in detail and, in particular, find
its  explicit structure. 

\thrml{281212m}
Let $\cA$ be a quasidense circulant S-ring.  Then $\cH(\cA)=\cH(\cA_{0,\frS})$ for any 
$\cA$-admissible class~$\frS$.
\ethrm
\proof Below to simplify notations we omit the letter $\frS$ in subscript. The theorem statement is obviously true when 
the S-ring $\cA$ is dense. Suppose that it is not dense. Then by Theorems~\ref{301012a} and~\ref{290113a} the S-ring 
$\cA$ admits an $\frS$-consistent $U/L$-decomposition. 

\lmml{230413a}
Let $\cB$ be a circulant S-ring over $G$. Suppose that $\cB$ is a $U/L$-wreath product. Then
$$
\cH(\cB)=\cH(\cB_U)\,\cup\,\pi^{-1}(\cH(\cB_{G/L}))
$$
where $\pi=\pi_{G/L}$ is the quotient epimorphism from $G$ to $G/L$.
\elmm
\proof Obviously, the right-hand side is contained in the left-hand  one. Conversely, let $H\in\cH(\cB)$. Without loss
of generality we can assume that $H\not\le U$. Then any highest basic set of $\cB_H$ is outside $U$.
Since $\cB$ is a $U/L$-wreath product, the radical of that set contains $L$. Thus $L\le H$ as required.\bull\medskip

By Lemma~\ref{230413a} and Lemma~\ref{241212a} with $S=G$
this implies that
$$
\cH(\cA)=\cH(\cA_U)\,\cup\,\pi^{-1}(\cH(\cA_{G/L}))\qaq \cH(\cA_0)=\cH((\cA_0)_U)\,\cup\,\pi^{-1}(\cH((\cA_0)_{G/L})).
$$
On the other hand, by induction we have 
\qtnl{080414a}
\cH(\cA_U)=\cH((\cA_U)_0)\qaq \cH(\cA_{G/L})=\cH((\cA_{G/L})_0).
\eqtn
Thus it suffices to verify that given $S\in\{U,G/L\}$ we have
\qtnl{230413c}
\cH((\cA_S)_0)=\cH((\cA_0)_S).
\eqtn
By Lemma~\ref{111113c} with $\cA=(\cA_S)_0$ and $\cA'=(\cA_0)_S$ all we need to prove is that
any elementary nondense section of~$(\cA_S)_0$ is an elementary nondense section of~$(\cA_0)_S$. 
However  by~\eqref{080414a}, any such section $T$ is an $\cA_S$-section. Therefore $T$
is an elementary nondense section of~$\cA$.  In particular, $\cA_T$ admits a unique $\frS$-consistent
gwr-decomposition, and by Theorem~\ref{290113a} this decomposition can be lifted to an $\frS$-consistent 
gwr-decomposition of~$\cA$.   By Lemma~\ref{241212a} this implies that $T$ is an elementary nondense section of~$\cA_0$, and
hence of $(\cA_0)_S$.\bull\medskip

It is easily seen that $(\cA_S)_{0,\frS}\le(\cA_{0,\frS})_S$ for all $\cA$-sections~$S$. The following
theorem refines this simple statement.

\thrml{291112a}
Let $\cA$ be a quasidense circulant S-ring and $\frS$ an $\cA$-admissible class. Then $(\cA_S)_{0,\frS}=(\cA_{0,\frS})_S$ for any
$\cA$-section~$S$.
\ethrm
\proof Below to simplify notations we omit the letter $\frS$ in subscript. Suppose on the contrary that 
$(\cA_S)_0<(\cA_0)_S$ for some $\cA$-section~$S$. Then there exist basic sets $X$ and $Y$ of the
S-rings $(\cA_S)_0$ and $(\cA_0)_S$ respectively, such that $Y$ is a proper subset of $X$.  From
Theorem~\ref{281212m} it follows that 
$$
\cH((\cA_0)_S)=\cH(\cA_S)=\cH((\cA_S)_0).
$$
Therefore $\grp{Y}$ is an $(\cA_S)_0$-group. However, the set $X\cap\grp{Y}$ is not empty. Thus,
$X\subset\grp{Y}$. On the other hand, $X$ and $Y$ are cosets,
because the S-rings $(\cA_S)_0$ and $(\cA_0)_S$ are coset ones. Since $Y$ is a proper subset of $X$,
this implies that  $\rad(Y)<\rad(X)$. Thus, $((\cA_S)_0)_T\ne\mZ T$
where $T=\grp{Y}/\rad(Y)$. This implies that
$$
(\cA_T)_0=((\cA_S)_T)_0\le ((\cA_S)_0)_T<\mZ T.
$$
By implication $(1)\Rightarrow(4)$ of Theorem~\ref{291112u} the hypothesis of Lemma~\ref{241212a}
is satisfied for $S=T$. Therefore the S-ring $(\cA_0)_T$ admits a gwr-decomposition.
It follows that $\rad((\cA_0)_T)\ne 1$, which  is impossible because $T$
is a principal section of the S-ring $(\cA_0)_S$.\bull

\crllrl{221113a}
In the conditions of Theorem~\ref{291112a} the following statements hold:
\nmrt
\tm{1} $\frSH=\frS_0(\cA_{0,\frS})$,
\tm{2} $\frSH$ is an $\cA$-admissible class  if and only if $\frSH\subset\frS_{cyc}(\cA)$,
\tm{3} $\frS_0=\frS_0(\cA_0)=\wh{\frS_0(\cA)}$ where $\frS_0$ and $\cA_0$ are as in Introduction.
\enmrt
\ecrllr
\proof Statement~(1) immediately follows from  equivalence $(1)\Leftrightarrow(2)$ of Theorem~\ref{291112u} 
and Theorems~\ref{291112a} and~\ref{270612b}. Statement~(2) follows from statement~(1) and
Theorem~\ref{281212m}. Statement~(3) is a special case of statement~(1) for $\frS=\frS_0(\cA)$.\bull\medskip

Let $\cA$ be a quasidense circulant S-ring and $\frS$ an $\cA$-admissible class. For a
basic set  $X$ of $\cA$ set
\qtnl{131113a}
L_\frS(X)=\bigcap_{\grp{X}/L\in\frSH, L\le\rad(X)} L.
\eqtn
Certainly, at least one group $L$ does exist because $\frSH\supset\frS_{prin}(\cA)$, and hence
one can take $L=\rad(X)$.
It should be mentioned that since the class $\frSH$ is closed with respect to taking subsections, the left-hand side 
of~\eqref{131113a} does not change when the intersection is taken over all 
sections $U/L\in\frSH$ having $\grp{X}/\rad(X)$ as a subsection.  Clearly, $L_\frS(X)$ is an $\cA$-group
contained in~$\rad(X)$. Therefore $X$ is a union of cosets of it; the set of all of them is 
denoted by~$X/L_\frS(X)$.

\thrml{281212b}
Let $\cA$ be a quasidense circulant S-ring and $\frS$ an $\cA$-admis\-sible class. Then
\nmrt
\tm{1} $\cS(\cA_{0,\frS})=\bigcup_{S\in\cS(\cA)}X/L_\frS(X)=\{xL_\frS(X):\ x\in X\in\cS(\cA)\}$,
\tm{2} elements $x$ and $y$ of a basic set $X$ of~$\cA$ are in the same basic set of $\cA_{0,\frS}$ if and only if
$\pi_S(x)=\pi_S(y)$ for any section $S\in\frSH$ having $\grp{X}/\rad(X)$ as a subsection.
\enmrt
\ethrm
\proof Statement~(2) immediately follows from statement~(1) and the remeark before the theorem. To prove statement~(1) 
let $X\in\cS(\cA)$ and $x\in X$. Denote by $X_0$ the basic set of the S-ring $\cA_{0,\frS}$ that contains $x$.
Then $X_0=xL_0$ for some $\cA_{0,\frS}$-subgroup~$L_0$, because $\cA_{0,\frS}$ is a coset  S-ring.
Moreover, by statement~(3) of Lemma~\ref {031212a} the set $X_0$ is contained in some $L$-coset in $U:=\grp{X}$
for any group $L$ such that $U/L\in\frSH$. By~\eqref{131113a} this implies that
$$
L_0\le L_\frS(X). 
$$

If $L_0=L_\frS(X)$, then $xL_\frS(X)=X_0$ and we are done.
Suppose that $L_0<L_\frS(X)$. By Theorem~\ref{281212m} the group $U_0:=\grp{X_0}$ is
an $\cA$-group. Since $X$ intersects $U_0$, this implies that $\grp{X}=U_0$. Moreover, due to the assumption 
we also have $L_0< \rad(X)$.
On the other hand, since $S:=U_0/L_0$ is a principal section of a coset S-ring $\cA_{0,\frS}$,
Theorem~\ref{270612b} implies that  $(\cA_{0,\frS})_S=\mZ S$. By 
Theorem~\ref{291112a} and implication (2)$\Rightarrow$(1) of Theorem~\ref{291112u}  this implies that
$S\in\frSH$. Since $\grp{X}/\rad(X)$ is a subsection of $S$, we obtain that  $L_\frS(X)\le L_0$. 
Contradiction.\bull

\section{Multipliers}\label{200214c}

Let $\cA$ be a quasidense circulant S-ring and $\frS$ an $\cA$-admissibe class.  In what
follows for an $\cA$-section $S$ we set $\aut_\cA(S)=\aut(\cA_S)\cap\aut(S)$.

\dfntnl{200313b}
An element $\fS=\{\sigma_S\}$ of the direct product $\prod_{S\in\frS}\aut_\cA(S)$, 
is called an $\frS$-multiplier of~$\cA$ if the following two conditions are satisfied for all sections $S_1,S_2\in\frS$:
\nmrt
\tm{M1} if $S_1\succeq S_2$, then $(\sigma_{S_1})^{S_2}=\sigma_{S_2}^{}$,
\tm{M2} if $S_1\sim S_2$ implies $m(\sigma_{S_1})=m(\sigma_{S_2})$.
\enmrt
The group of all $\frS$-multipliers of~$\cA$ is denoted by $\mult_\frS(\cA)$.
\edfntn

It should be noted that if the class  $\frS_0$ defined in~\eqref{110214b} is contained in $\frS_{cyc}(\cA)$, then 
it is admissible (Corollary~\ref{221113a}) and $\mult(\cA)=\mult_{\frS_0}(\cA)$. The following lemma gives
a natural way to construct $\frS$-multipliers.

\lmml{140414a}
Suppose that we are given $\gamma\in\aut(\cA)$ such that $\gamma^S\in\aut_\cA(S)$
for all $S\in\frS$.
Then the family $\fS(\gamma)=\{\sigma_S(\gamma)\}_{S\in\frS}$ where $\sigma_S(\gamma)=\gamma^S$,
is an $\frS$-multiplier.
\elmm
\proof Condition~(M1) is obvious. Besides, from Theorem~\ref{261010a} 
it follows that $m(\gamma^{S_1})=m(\gamma^{S_2})$ for any projectively equivalent
sections $S_1,S_2\in\frS$. Therefore $m(\sigma_{S_1}(\gamma))=m(\sigma_{S_2}(\gamma))$. Thus condition
(M2) is also satisfied for $\fS(\gamma)$, and we are done.\bull\medskip

Let $\fS\in\mult_\frS(\cA)$ and $T\in\frS(\cA)$. Then $\frS_T$
is an $\cA_T$-addmissible class by statement~(2) of Lemma~\ref{141113v}, and
$\frS_T\subset\frS$ by~\eqref{281113a}. Conditions (M1) and (M2) 
are obviously satisfied for the restriction $\fS^T$ of $\fS$ to $\frS_T$, which is by definition
considered as an element of the direct product $\prod_{S\in\frS_T}\aut_\cA(S)$. This
proves the following statement.


\lmml{261113b}
The family $\fS^T$ is an $\frS_T$-multiplier of the S-ring $\cA_T$.\bull
\elmm

We are going to construct $\frS$-multipliers of $\cA$ by means of similarities belonging to the set
$$
\Phi_{0.\frS}(\cA)=\{\varphi\in\Phi(\cA_{0,\frS}):\ X^\varphi=X\ \text{for all}\ X\in\cS(\cA)\}.
$$
Namely, with any $\varphi\in \Phi_{0,\frS}(\cA)$ we associate an element $\fS_\varphi=\{\sigma_S\}$
of the group $\prod_{S\in\frS}\aut_\cA(S)$ where the automorphisms $\sigma_S$ are
defined as follows. Let $S\in\frS$.  Then by statement~(2) of Lemma~\ref{031212a}
we have $(\cA_{0,\frS})_S=\mZ S$. Therefore   there exists a uniquely determined automorphism 
$\sigma_S\in\aut(S)$ that induces the restriction $\varphi_S$ of the similarity~$\varphi$ to~$S$  
(Lemma~\ref{271213a}). From the choice of $\varphi$ it follows that $\sigma_S\in\aut_\cA(S)$. 

\lmml{281113d}
The family $\fS_\varphi$ is an $\frS$-multiplier of the S-ring $\cA$. Moreover,
$\varphi_T\in\Phi_{0,\frS_T}(\cA_T)$ and $(\fS_\varphi)^T=\fS_{\varphi_{\scriptscriptstyle T}}$ for all $T\in\frS(\cA)$.
\elmm
\proof By Theorem~\ref{150113a} the similarity $\varphi$ is induced by an isomorphism $\gamma$ of the S-ring
$\cA_{0,\frS}$. Clearly, $\gamma$ can be chosen so that $1^\gamma=1$.
Then by Lemma~\ref{271213a} we have  $\gamma^S\in\aut_\cA(S)$ for all $S\in\frS$. Besides, 
$\gamma\in\aut(\cA)$ by Lemma~\ref{021213a}. So 
from Lemma~\ref{140414a} it follows that $\fS(\gamma)$  is an $\frS$-multiplier 
of~$\cA$. Thus, the first statement follows because $\fS_\varphi=\fS(\gamma)$.\medskip


To prove the second statement, let $T\in\frS(\cA)$. Then $\varphi_T\in \Phi((\cA_{0,\frS})_T)$. 
By Theorem~\ref{291112a}
 this implies that $\varphi_T\in\Phi((\cA_T)_{0,\frS})$. So $\varphi_T\in\Phi_{0,\frS_T}(\cA_T)$.
 The rest of the statement easily follows from the definitions.\bull\medskip

 To simplify notation we will write $\Phi_{0,\frS}(\cA_T)$ instead of $\Phi_{0,\frS_T}(\cA_T)$.
 
\thrml{071212d}
Let $\cA$ be an quasidense circulant S-ring. Then the mapping
\qtnl{221113w}
\Phi_{0,\frS}(\cA)\to \mult_\frS(\cA),\quad \varphi\mapsto\fS_\varphi
\eqtn
is a group isomorphism for any $\cA$-admissible class~$\frS$.
\ethrm
\proof The mapping~\eqref{221113w} is obviously a group homomorphism. To prove its injectivity
suppose that $\fS_\varphi=\fS_\psi$ for some $\varphi,\psi\in\Phi_{0,\frS}(\cA)$. Then obviously
$\varphi_S=\psi_S$ for all $S\in\frS$. By Lemma~\ref{271213a} the equality also holds for all $S\in\frSH$ 
because any similarity of $\mZ S$ is uniquely determined by its restrictions to the Sylow sections $S_p$.
However, $\frSH=\frS_0(\cA_{0,\frS})$  by Corollary~\ref{221113a}. 
Therefore $\varphi$ and $\psi$ are equal on
all principal $\cA_{0,\frS}$-sections. Thus $\varphi=\psi$ by Lemma~\ref{221113f}.\medskip

Let us prove the surjectivity of homomorphism~\eqref{221113w} by induction on the size of the group~$G$ 
underlying~$\cA$. Let $\fS=\{\sigma_S\}$ be an $\frS$-multiplier of~$\cA$. First, suppose that
$\cA_{0,\frS}=\mZ G$. By  implication $(1)\Rightarrow(2)$ of Theorem~\ref{291112u}
this implies that $G\in\frSH$. So any Sylow subgroup of $G$ belongs to~$\frS$. It follows that $\cA$ is a subtensor 
product of cyclotomic S-rings. In this case the following statement holds.

\lmml{071212e}
There exists  $\gamma\in\aut_\cA(G)$ such that $\fS=\fS(\gamma)$.
\elmm
\proof Set $\gamma$ to be the unique automorphism of $G$ such that $\gamma^{G_p}=\sigma_{G_p}$ for all Sylow subgroups~$G_p$ of~$G$.
We claim that
\qtnl{071212f}
\sigma_S(\gamma)=\sigma_S,\quad S\in\frS,
\eqtn
where $\sigma_S(\gamma)$ is the $S$-component of the family $\fS(\gamma)$.
Indeed, the automorphism $\sigma_S$ is uniquely determined by its $p$-components $(\sigma_S)_p\in\aut(S_p)$. 
Besides, $S_p\in\frS$ because $S\in\frS$, and hence $(\sigma_S)_p=\sigma_{S_p}$ by  Definition~\ref{200313b}.
Since $G_p\in\frS$,  each $\sigma_{S_p}$ is in its turn uniquely determined by the automorphism $\sigma_{G_p}$.
By the definition of~$\gamma$ this proves~\eqref{071212f}. \medskip

To complete the proof of the lemma let us verify that $\gamma\in\aut(\cA)$.
Since $\gamma\in\aut(G)$,  by Lemma~\ref{021213a} for $\cA'=\mZ G$, it suffices to check that 
$X^\gamma=X$ for all basic sets~$X$ of~$\cA$.
However, the section $S=\grp{X}/\rad(X)$ belongs to~$\frS$. Moreover, $X$ is a disjoint union of $\rad(X)$-cosets.
Since $\sigma_S\in\aut(\cA_S)$, equality~\eqref{071212f} implies that the automorphism $\gamma$ permutes  
the cosets in this union, as required.\bull\medskip

By Lemma~\ref{071212e} we have $\fS=\fS(\gamma)$. On the other hand, it is easily seen that
$\fS(\gamma)=\fS_\varphi$ where $\varphi=\varphi_\gamma$  is the
similarity induced by the automorphism~$\gamma$ (Lemma~\ref{271213a}). Thus
$\fS=\fS_\varphi$  as required.\medskip

Now, assume that $\cA_{0,\frS}\ne\mZ G$.
Then by implication (1)$\Rightarrow$(3) of Theorem~\ref{291112u} for $S=G$.
there exists an $\frS$-consistent  $U/L$-decomposition of $\cA$.
By the inductive hypothesis applied to the S-ring $\cA_U$ and $\frS_U$-multiplier $\fS^U$, as well
to  the S-ring $\cA_{G/L}$ and $\frS^{G/L}$-multiplier $\fS^{G/L}$, there exist similarities
$\varphi_1\in\Phi_{0,\frS}(\cA_U)$ and $\varphi_2\in\Phi_{0,\frS}(\cA_{G/L})$
such that 
\qtnl{291113b}
\fS^U=\fS_{\varphi_1}\qaq\fS^{G/L}=\fS_{\varphi_2}.
\eqtn
 Furthermore, $(\cA_U)_{0,\frS}=(\cA_{0,\frS})_U$ and $(\cA_{G/L})_{0,\frS}=(\cA_{0,\frS})_{G/L}$ by Theorem~\ref{291112a}.
Therefore
\qtnl{110414a}
((\cA_U)_{0,\frS})_{U/L}=(\cA_{0,\frS})_{U/L}=((\cA_{G/L})_{0,\frS})_{U/L}.
\eqtn
However, again by Theorem~\ref{291112a} we have $(\cA_{0,\frS})_{U/L}=(\cA_{U/L})_{0,\frS}$. It follows that
the similarities $(\varphi_1)^{U/L}$ and $(\varphi_2)^{U/L}$ belong to $\Phi_{0,\frS}(\cA_{U/L})$.\medskip

By the first equality in~\eqref{291113b} and Lemma~\ref{281113d} applied to $\cA=\cA_U$, $\varphi=\varphi_1$ and $T=U/L$,
we obtain that
$$
\fS_{(\varphi_1)^{U/L}}=(\fS_{\varphi_1})^{U/L}=(\fS^U)^{U/L}=\fS^{U/L}.
$$
Similarly, by the second equality in~\eqref{291113b} and Lemma~\ref{281113d} applied to $\cA=\cA_U$, $\varphi=\varphi_1$ 
and $T=U/L$ we have
$$
\fS_{(\varphi_2)^{U/L}}=(\fS_{\varphi_2})^{U/L}=(\fS^{G/L})^{U/L}=\fS^{U/L}.
$$
Thus $\fS_{(\varphi_1)^{U/L}}=\fS_{(\varphi_2)^{U/L}}$ and by the injectivity statement we have
$$
(\varphi_1)^{U/L}=(\varphi_2)^{U/L}.
$$
On the other hand, by Theorem~\ref{230114a} and equality~\eqref{110414a} one can form the 
$U/L$-wreath product $\cA'$ of the S-rings $(\cA_U)_{0,\frS}$ and 
$(\cA_{G/L})_{0,\frS}$. Thus by  statement~(2) of Theorem~\ref{271213x} there exists 
a uniquely determined similarity $\varphi\in\Phi(\cA')$ such that 
\qtnl{110414d}
\varphi^U=\varphi_1\qaq\varphi^{G/L}=\varphi_2.
\eqtn
However, $\cA$ is the $U/L$-wreath product. By Corollary~\ref{241212a} so is $\cA_{0,\frS}$. Since the
restrictions of the latter S-ring to $U$ and $G/L$ coincide with the corresponding restrictions of the
S-ring~$\cA'$ (see above), we conclude that $\cA'=\cA_{0,\frS}$. Thus $\varphi\in\Phi(\cA_{0,\frS})$. Since
$\varphi_1\in\Phi_{0,\frS}(\cA_U)$ and $\varphi_2\in\Phi_{0,\frS}(\cA_{G/L})$, we conclude that
$\varphi\in\Phi_{0,\frS}(\cA)$.\medskip

To complete the proof let us verify that $\fS=\fS_\varphi$. To do this we observe that from~\eqref{291113b}
and~\eqref{110414d} it follows that $(\fS_\varphi)^U=\fS^U$ and $(\fS_\varphi)^{G/L}=\fS^{G/L}$.
This proves the required statement because by the $\frS$-consistency property, any section in~$\frS$
is either $\cA_U$- or $\cA_{G/L}$-section.\bull

\section{Proof of Theorems~\ref{251212a} and~\ref{170113a}}\label{200214f}

The following auxiliary statement is interesting by itself. In particular, it shows that the subgroup lattices
of a quasidense S-ring and its schurian closure are equal. It seems that this statement could be generalized
to all circulant  S-rings.

\lmml{230415e}
Let $\cA$ be a quasidence cirulant S-ring and $\cA'$ its schurian closure. Then $\cH(\cA)=\cH(\cA')$ and 
$\cA_0=(\cA')_0$.
\elmm
\proof From Theorem~\ref{150113a} it follows that any coset S-ring that contains $\cA$, contains also $\cA'$.
This proves the second equality and shows that  $\cA_0\ge\cA'$. Since 
obviously $\cA'\ge\cA$, the first equality follows from equality~\eqref{191113s} and Theorem~\ref{281212m}.\bull\medskip

{\bf Proof of Theorem~\ref{251212a}.} Without loss of generality we can assume that the S-ring $\cA$ is schurian. 
Indeed, from Lemma~\ref{230415e} it follows that $\cA_0=(\cA')_0$ where $\cA'=\sch(\cA)$. Therefore it suffices to verify 
that
\qtnl{021213d}
\Phi_0(\cA)=\Phi_0(\cA').
\eqtn
 However, by Theorem~\ref{150113a} any similarity $\varphi\in\Phi_0(\cA)$ is induced by an isomorphism 
$f$ of $\cA_0$ to itself. Without loss of generality we can assume that $1^f=1$. Then by~\eqref{260613a} 
and the choice of~$\varphi$
we have $X^f=X^\varphi=X$ for all $X\in\cS(\cA)$. By Lemma~\ref{021213a} this implies that 
$f\in\aut(\cA)$, and hence $f\in\aut(\cA')$. Then $\varphi$ is
 identical on~$\cA'$, and so $\varphi\in\Phi_0(\cA')$. 
Thus $\Phi_0(\cA)\subset\Phi_0(\cA')$. Since the reverse inclusion is obvious,
equality~\eqref{021213d} follows.\medskip

From the definition of the group $\Phi_0=\Phi_0(\cA)$ it follows that $\cA\le \cB$ where $\cB=(\cA_0)^{\Phi_0}$. To prove 
the reverse inclusion we have to verify that any basic set of $\cA$ is contained in a basic set of $\cB$. Let
$x$ and $y$ belong to the same basic set of $\cA$. Then it suffices to verify that 
\qtnl{021213e}
x^\delta=y\quad\text{for some}\quad \delta\in\aut(\cB),\ \ 1^\delta=1.
\eqtn
To do this we recall that $\cA$ is schurian. So by Theorem~\ref{230312a} one can find a group $\Gamma\in\cM(\cA)$ 
such that
\qtnl{150113b}
\Gamma^S=\hol_\cA(S),\qquad S\in\frS_0(\cA).
\eqtn
It follows that there exists $\gamma\in\Gamma$ such that $1^\gamma=1$ and $x^\gamma=y$. Due to
\eqref{150113b} we have $\gamma^S\in\aut_\cA(S)$ for all $S\in\frS_0(\cA)$. 
By Lemma~\ref{140414a} this implies that the family 
$\fS=\{\gamma^S\}$ is an $\frS_0(\cA)$-multiplier of~$\cA$.
Therefore by Theorem~\ref{071212d} there exists  a uniquely determined similarity $\varphi\in\Phi_0$ such that
$\fS=\fS_\varphi$. This means that $\varphi^S$ is induced by  $\gamma^S$ for all $S\in\frS_0(\cA)$. Since 
the S-ring $\cA_0$ is separable,
there exists an isomorphism $\gamma_0$ of the ring~$\cA_0$ to itself that induces~$\varphi$. Without loss of 
generality we can assume that $1^{\gamma_0}=1$. We claim that
\qtnl{021213f}
x^{\gamma_0}=y^{\gamma'}\quad\text{for some}\quad \gamma'\in\aut(\cA_0),\ \ 1^{\gamma'}=1.
\eqtn
Then \eqref{021213e} holds for $\delta=\gamma_0(\gamma')^{-1}$ because $\gamma_0\in\aut(\cB)$ by 
Lemma~\ref{021213a}, and we are done.\medskip

Let us prove the claim. From Lemma~\ref{021213a} it follows that $\gamma_0\in\aut(\cB)$. Therefore $x^{\gamma_0}\in X$ where $X$
is the basic set of~$\cA$ that contains $x$ and $y$. Let $S\in\frS_0(\cA)$. Then the bijections $\gamma_0^S$
and $\gamma^S$ induce the same similarity $\varphi^S$ of the S-ring $(\cA_0)_S$. However, this S-ring equals 
$\mZ S$ by statement~(2) of Lemma~\ref{031212a}. Therefore $\gamma_0^S=\gamma^S$. 
Thus the latter equality holds for all $S\in\frS_0(\cA)$.\medskip

Let now $S$ be a section  in the class $\wh{\frS_0(\cA)}$ that has $\grp{X}/\rad(X)$ as a subsection.
Then $S_p\in\frS_0(\cA)$ for all primes~$p$ dividing $|S|$. By above this implies that
$$
x_p^{\gamma_0}L=(x_pL)^{\gamma_0}=(x_pL)^\gamma=x_p^\gamma L=y_pL
$$
where $L$ is the denominator of~$S$.
On the other hand, the S-ring $\cA_S$ contains the tensor product of the S-rings 
$\cA_{S_p}$, and hence the permutation $(\gamma_0)^S$ equals the product of its $p$-projections. Thus 
$x^{\gamma_0}L=yL$. By statement~(2) of Theorem~\ref{281212b} this implies that the elements 
$x^{\gamma_0}$ and $y$ belong to  the same basic set of $\cA_0$. Since this S-ring is schurian, there
exists an automorphism $\gamma'\in\aut(\cA_0)$ such that $x^{\gamma_0}=y^{\gamma'}$ and $1^{\gamma'}=1$.
The claim is proved.\bull\medskip

{\bf Proof of Theorem~\ref{170113a}.} To prove the ``only if'' part suppose that the S-ring $\cA$ is schurian. Then by 
Theorem~\ref{230312a} there exists a group $\Gamma\in\cM(\cA)$ that  
satisfies~\eqref{150113b}. Let $S\in\frS_0$. Then  since $\frS_0=\wh{\frS_0(\cA)}$
(Corollary~\ref{221113a}), we have $S_p\in\frS_0(\cA)$ for any prime $p$ dividing $|S|$. Therefore due to ~\eqref{150113b} we conclude that
$$
\Gamma^S\le\prod_p\Gamma^{S_p}\le\prod_p\hol(S_p)=\hol(S).
$$
However, $\cA_S$ is the S-ring associated with the group $\Gamma^S$ because $\cA$ is the S-ring associated with 
the group $\Gamma$.  Thus the S-ring $\cA_S$ is cyclotomic and condition~(1) is satisfied. To verify condition~(2) we 
note that $\aut_\cA(S)=\Gamma^S$. Therefore given $\sigma\in\aut_\cA(S)$ one can find $\gamma\in\Gamma$ such that $1^\gamma=1$ and 
$\gamma^S=\sigma$. Thus the element $\sigma$ has a preimage  $\fS=\{\gamma^S\}$ in the 
group~$\mult(\cA)$.\medskip

To prove the "if part" suppose that  for any section $S\in\frS_0$  the S-ring $\cA_S$ is cyclotomic and
the restriction homomorphism from $\mult(\cA)$ to $\aut_\cA(S)$ is surjective. Then by the first assumption
the class $\frS_0$ is an $\cA$-admissible (Corollary~\ref{221113a}) and 
\qtnl{180113a}
\cA_S=\cyc(\aut_\cA(S),S)=(\mZ S)^{\aut_\cA(S)}
\eqtn
for all $S\in\frS_0$. Moreover, we claim that
\qtnl{260413b}
(\mZ S)^{\aut_\cA(S)}=(\mZ S)^{(\Phi_0)^S}
\eqtn
where $\Phi_0=\Phi_0(\cA)$.
Indeed, the left-hand side is obviously contained in the right-hand side. Suppose on the contrary that the reverse 
inclusion does not hold. Then there exists $\sigma\in \aut_\cA(S)$ such that $\varphi_\sigma\not\in (\Phi_0)^S$ 
where $\varphi_\sigma$ is the similarity of $\mZ S$ induced by~$\sigma$. Moreover, by the surjectivity assumption
there is a family $\Sigma\in \mult(\cA)$  the $S$-component of which equals~$\sigma$. 
Since the class $\frS_0$ 
is an $\cA$-admissible,  $\Sigma$ is an $\frS_0$-mulliplier of $\cA$. So by Theorem~\ref{071212d} there exists a similarity $\varphi\in\Phi_0$ such that
$\varphi_S=\varphi_\sigma$. This implies that  $\varphi_\sigma\in (\Phi_0)^S$. Contradition.\medskip

Set $\cA'=(\cA_0)^{\Phi_0}$. Then from~\eqref{180113a} and \eqref{260413b} it follows that
\qtnl{260413a}
(\cA')_S=((\cA_0)^{\Phi_0})_S=((\cA_0)_S)^{(\Phi_0)^S}=(\mZ S)^{(\Phi_0)^S}=\cA_S
\eqtn
for all $S\in\frS_0$. By Theorem~\ref{251212a} to complete the proof it suffices to verify that $\cA=\cA'$. To do this 
we note
that $\cA\le \cA' $.  So we only need to prove that $\cA\ge\cA'$, i.e. that any basic set $X'\in\cS(\cA')$ belongs 
to $\cS(\cA)$. However, obviously $S'=\grp{X'}/\rad(X')$ is an $\cA_0$-section. Moreover, the extension
$(\cA_0)_{S'}$ of $(\cA')_{S'}$ has trivial radical by Theorem~\ref{060214b}. By Lemma~\ref{100214a}
this implies that
$S'\in\frS_0$. Denote by $X$ the basic set of the S-ring $\cA_{\grp{X'}}$ that contains~$X'$. Then since
$\cA^{}_{S'}=\cA'_{S'}$ (see~\eqref{260413a}), we have $\pi_{S'}(X)=\pi_{S'}(X')$. It follows that
 $X\subset X'\rad(X')=X'$. Thus $X=X'$  as required.\bull

\section{Proof of Theorem~\ref{170113b}}\label{100214b}

\sbsnt{Reduction to quasidense S-rings.}
Let $\cA$ be an S-ring over a cyclic group $G$. Following paper~\cite{EP12} a class  $C$ of projectively equivalent 
$\cA$-sections is called {\it singular} if its  rank is~$2$, its order is greater than~$2$ and it contains two 
sections $S=L_1/L_0$ and $T=U_1/U_0$ such that $T$ is a multiple of~$S$ and
the following two conditions are satisfied:
\nmrt
\tm{S1} $\cA$ is both the $U_0/L_0$- and $U_1/L_1$-wreath product,
\tm{S2} $\cA_{U_1/L_0}=\cA_{L_1/L_0}\otimes\cA_{U_0/L_0}$.
\enmrt
By Lemma~6.2 of that paper $S$ and  $T$ are necessarily  the smallest and largest $\cA$-sections of $C$. 
Moreover, from Theorem~4.6 there, it follows 
that any rank~$2$ class of  composite order belonging to $\cP(\cA)$ is singular; in particular,  the S-ring $\cA$ is
quasidense if and only if no class in $\cP(\cA)$ is singular.\medskip 

For an $\cA$-section $S$ we define the {\it $S$-extension} of~$\cA$ to be the  smallest 
S-ring $\cA'\ge\cA$ such that $\cA'_S=\mZ S$. From Theorem~\ref{261010a} it follows that $\cA'$
does not depend on the choice of $S$ in the class $C\in\cP(\cA)$ of sections projectively equivalent to~$S$.
Denote this S-ring by $E(\cA,C)$. \medskip

 The following lemma provides the reduction of Theorem~\ref{170113b}
to the quasidense case, and will be used later only for a singular class~$C$ of composite order.

\lmml{101213a}
Let $\cA$ be a circulant S-ring, $C\in\cP(\cA)$ a singular class and $\cA'=E(\cA,C)$. Then
\nmrt
\tm{1} $\rk(\cA')>\rk(\cA)$,
\tm{2} $\cA'$ is schurian if and only if $\cA$ is schurian.
\enmrt
\elmm
\proof Statement~(1) follows from the fact that any singular class has rank~2 and order at least~3.
To prove statement~(2) denote by $L_1/L_0$ and  $U_1/U_0$ the smallest and the largest $\cA$-sections of $C$
(see Theorem~\ref{040214b}). 
First, we claim that  $\cA'$ coincides with any extension $\cB$ of $\cA$ that satisfies
the following conditions 
\nmrt
\tm{E1} $\cB$ is both the $U_0/L_0$- and $U_1/L_1$-wreath product,
\tm{E2} $\cB_{U_0}=\cA^{}_{U_0}$, $\cB_{G\setminus L_1}=\cA^{}_{G\setminus L_1}$ and 
$\cB_{U_1/L_0}=\mZ S\otimes\cA^{}_{U_0/L_0}$
\enmrt
where $S=L_1/L_0$. Indeed, from condition~(E2) it follows that $\cB_S=\mZ S$.
So by the definition of $\cA'$ we have
$$
\cA\le\cA'\le\cB\qaq (\cA')_S=\mZ S=\cB_S.
$$
This implies that conditions (E1) and (E2) are satisfied with $\cB$ replaced by~$\cA'$. It follows that
\qtnl{160414a}
\cB^{}_{U_0}=\cA'_{U_0},\quad \cB^{}_{G\setminus L_1}=\cA'_{G\setminus L_1},\quad
\cB^{}_{U_1/L_0}=\cA'_{U_1/L_0},
\eqtn
and $\cA'$ and $\cB$ are the $U_1/L_1$- and $U_0/L_0$-wreath products. Therefore
to check that $\cA'=\cB$ it suffices to verify that 
$$
\cA'_{G/L_1}=\cB^{}_{G/L_1}\qaq
\cA'_{U_1}=\cB^{}_{U_1}.
$$
The former equality immediately follows from the second equality in~\eqref{160414a}. 
To verify the latter one, we observe that $\cA'_{U_1}$ and $\cB^{}_{U_1}$ are the $U_0/L_0$-wreath products. 
Thus the required statement follows from the other two equalities in~\eqref{160414a}. \medskip

To complete the proof we note that in terms of paper~\cite{EP12} every singular class is 
isolated (i.e. satisfies conditions~(S1) and~(S2)). Therefore by Lemma~6.5 of that paper the S-ring $\Ext_C(\cA,\mZ S)$ defined there
contains $\cA$ and satisfies conditions~(E1) and~(E2). By the above claim this implies
that $\Ext_C(\cA,\mZ S)=\cA'$. Thus  statement~(2) follows from \cite[Theorem~6.7]{EP12}.\bull\medskip

We recall that $\wh\cA$  is the S-ring dual to~$\cA$, and $\wh C$ is the class of projectively equivalent $\wh\cA$-sections 
that is dual to a class $C\in\cP(\cA)$  (Subsection~\ref{091213a}).

\thrml{270612a}
Let $\cA$ be a circulant S-ring and $C$ a singular class of $\cA$. Then $\wh C$ is a singular class 
of $\wh\cA$ and the S-ring dual to $E(\cA,C)$ coincides with $E(\wh\cA,\wh C)$.
\ethrm
\proof  Let $S\in C$ be a section of rank~$2$  and order at least~$3$. Then the class~$\wh C$ contains 
the section $\wh S$. Since $|S|=|\wh S|$ and $\rk(\wh\cA_{\wh S})=\rk(\cA_S)=2$
(Subsection~\ref{091213a}), the class $\wh C$ has rank~$2$ and order greater than~$2$. Moreover,
$U_0^\bot/U_1^\bot$ and $L_0^\bot/L_1^\bot$ are $\wh\cA$-sections, and statement~(1) of Lemma~\ref{310114a} implies that 
the former one is a multiple of the latter. Finally, by Theorem~\ref{140509a} the S-ring~$\wh\cA$
satisfies conditions~(S1) and~(S2). Thus the class $\wh C$ is singular.\medskip

Let us prove that
\qtnl{101213c}
E(\wh\cA,\wh C)\le\wh{E(\cA,C)}.
\eqtn
Indeed, since $E(\cA,C)\ge\cA$, the S-ring dual to $E(\cA,C)$ contains $\wh\cA$. Moreover,
since $E(\cA,C)_S=\mZ S$, the restriction of that ring to $S^\bot$ equals $\mZ S^\bot$.
Thus~\eqref{101213c} follows from the definition of $E(\wh\cA,\wh C)$. Next by duality,
inclusion~\eqref{101213c}  also holds after interchanging $\cA$ and $\wh\cA$. Therefore
$$
\wh{E(\cA,C)}\le E(\wh\cA,\wh C).
$$
Due to inclusion~\eqref{101213c} this completes the proof of the theorem.\bull\medskip


Let us turn to the proof of Theorem~\ref{170113b}. Let $\cA$ be a circulant S-ring. First,
we observe that  given a singular class 
$C\in\cP(\cA)$ of composite order
the S-ring $E(\cA,C)$ is schurian if and only if so is $\cA$ (statement~(2)
of Lemma~\ref{101213a}). Therefore by statement~(1) of that lemma and by
Theorem~\ref{270612a} without loss of generality we can assume that the S-ring $\cA$ has no singular classes
of composite order, or equivalently that $\cA$ is a quasidense S-ring (Theorem~\ref{060214a}).\medskip

\sbsnt{Quasidense S-ring case.}
Given $\sigma\in\aut(G)$ we set $\wh\sigma$ to be the automorphism of $\wh G$  taking a character $\chi$ to
the character  $\chi^{\wh\sigma}:g\mapsto\chi(g^\sigma)$.  Then $m(\sigma)=m(\wh\sigma)$ because
$$
\chi^{\wh\sigma}(g)=\chi(g^\sigma)=\chi(g^m)=\chi(g)^m,\qquad g\in G,
$$
where $m=m(\sigma)$. Moreover, this shows that $\wh{\sigma^S}=\wh\sigma^{\wh S}$ for any section~$S$ of~$G$.

\lmml{280113a}
Let $\cA$ be a quasidense circular S-ring. Then
\nmrt
\tm{1} $\frS_0(\wh\cA)=\{\wh S:\ S\in\frS_0(\cA)\}$,
\tm{2} $\wh{\cA_0}=\wh\cA_{\,0}$ and $\frS_0(\wh\cA_{\,0})=\{\wh S:\ S\in\frS_0(\cA_0)\}$,
\tm{3} $\aut_{\wh\cA}(\wh S)=\{\wh\sigma:\ \sigma\in\aut_{\cA}(S)\}$ for all $S\in\frS(\cA)$.
\enmrt
\elmm
\proof  Statement~(1) follows from Theorem~\ref{101213e}. From this statement and equivalence $(1)\Leftrightarrow(4)$
of Theorem~\ref{270612b} it follows that a circulant S-ring is coset if so is its dual. Therefore the set of duals to 
coset S-rings containing~$\cA$
coincides with the set of coset S-rings containing~$\wh\cA$. Thus $\wh{\cA_0}=\wh\cA_{\,0}$  by the definition of
the coset closure. The second part of statement~(2) follows from the first one and statement~(1).\medskip

To prove statement~(3) let $S\in\frS(\cA)$. Then the group $\aut_{\cA}(S)$ equals the largest
group $K\le\aut(S)$ for which $\cyc(K,S)\ge\cA_S$. Since the dual to the S-ring $\cyc(K,S)$ equals
$\cyc(\wh K,\wh S)$ where $\wh K=\{\wh\sigma:\ \sigma\in K\}$ (Theorem~\ref{210114a}), the group $\wh K$ is the 
largest subgroup of $\aut(\wh S)$ for which  $\cyc(\wh K,\wh S)\ge\wh\cA_{\wh S}$. So $\aut_{\wh\cA}(\wh S)=\wh K$
and we are done.\bull\medskip

To complete the proof of Theorem~\ref{170113b} let $\cA$ be a  schurian quasidence circulant S-ring. By
duality we only have  to prove that the S-ring $\wh\cA$ is schurian. However, the latter ring
is quasidense by statement~(2) of Theorem~\ref{201213a}.  Therefore by Theorem~\ref{170113a} 
 and statement~(3) of Corollary~\ref{221113a} it suffices to 
verify that both of its conditions are satisfied for the S-ring~$\wh\cA$ and sections belonging
to $\frS_0(\wh\cA_{\,0})$.\medskip

Let $T\in\frS_0(\wh\cA_{\,0})$. Then by statement~(2) of Lemma~\ref{280113a} there exists a section 
$S\in\frS_0(\cA_0)$ such that $T=\wh S$. Since $\cA$ is schurian, the S-ring $\cA_S$ is cyclotomic 
(Theorem~\ref{170113a}). Therefore the S-ring $\wh\cA_T=\wh{\cA_S}$ is cyclotomic by Theorem~\ref{210114a}.
Thus condition~(1) is satisfied.\medskip
 
To verify condition~(2),
let $\tau\in\aut_{\wh\cA}(T)$. Since $T\in\frS(\wh\cA)$  (Theorem~\ref{281212m}),
statement~(3) of Lemma~\ref{280113a} implies that
there exists $\sigma\in\aut_\cA(S)$ such that $\tau=\wh{\sigma}$.  Moreover, since $\cA$ is schurian
there exists an element $\Sigma=\{\sigma_{S'}\}$ of the group $\mult(\cA)$ such that $\sigma_S=\sigma$
(Theorem~\ref{170113a}).
Set
$$
\wh\Sigma=\{\sigma_{T'}\}_{T'\in\frS_0(\wh\cA_{\,0})}
$$
where $\sigma_{T'}=\wh{\sigma_{S'}}$ with $S'$ defined by $\wh{S'}=T'$ (see statement~(2) of
Lemma~\ref{280113a}). Since $\wh{\sigma_S}=\wh\sigma=\tau$,  
it suffices to verify that $\wh\Sigma\in\mult(\wh\cA)$. 
By the remark after Definition~\ref{200313b} we have to verify that conditions~(M1) and~(M2) of
that definition are satisfied for the class $\frS_0(\wh\cA_{\,0})$ and family
$\wh\Sigma$.\medskip

Let
$\tau_1\in\aut_{\wh\cA}(T_1)$ where $T_1\in \frS_0(\wh\cA_{\,0})$. Then since $\Sigma$ is a multiplier, 
for every subsection $T_2$ of $T_1$ we have
$$
(\tau_1)^{T_2}=(\wh{\sigma_1})^{\wh{S_2}}=\wh{(\sigma_1)^{S_2}}=\wh{\sigma_2}=\tau_2
$$ 
 where $S_i$ is such that $T_i=\wh{S_i}$ and
$\sigma_i\in\aut_\cA(S_i)$ is such that $\tau_i=\wh{\sigma_i}$, $i=1,2$. 
Thus condition (M1) is satisfied. To verify condition~(M2)
let now $T_2\sim T_1$. Then by statement~(1) of Lemma~\ref{310114a} we have $S_1\sim S_2$. 
Since $\Sigma$ is a multiplier, this implies that
$m(\sigma_1)=m(\sigma_2)$.  Thus
$$
m(\tau_1)=m(\wh{\sigma_1})=m(\sigma_1)=m(\sigma_2)=m(\wh{\sigma_2})=m(\tau_2)
$$ 
as required.

\end{document}